\documentclass[12pt]{amsart}
\font\emailfont=cmtt10

\usepackage{amsmath,amsthm,amsfonts,amscd,flafter,epsf}

\newcommand\commentable[1]{#1}

\newcommand{\Tors}{\mathrm{Tors}}
\newcommand{\rk}{\mathrm{rk}}

\newtheorem{theorem}{Theorem}[section]
\newtheorem{prop}[theorem]{Proposition}
\newtheorem{cor}[theorem]{Corollary}

\newtheorem{lemma}[theorem]{Lemma}

\newtheorem{defn}[theorem]{Definition}

\newtheorem{remark}[theorem]{Remark}

\def\endproof{\relax\ifmmode\expandafter\endproofmath\else
  \unskip\nobreak\hfil\penalty50\hskip.75em\hbox{}\nobreak\hfil\bull
  {\parfillskip=0pt \finalhyphendemerits=0 \bigbreak}\fi}
\def\endproofmath$${\eqno\bull$$\bigbreak}
\def\bull{\vbox{\hrule\hbox{\vrule\kern3pt\vbox{\kern6pt}\kern3pt\vrule}\hrule}}

\newcommand{\Q}{\mathbb{Q}}
\newcommand{\R}{\mathbb{R}}

\newcommand{\Z}{\mathbb{Z}}

\newcommand{\Zmod}[1]{\Z/{#1}\Z}

\newcommand{\cm}{\cdot}

\newcommand{\ModSWfour}{\mathcal{M}}
\newcommand{\ModFlow}{\ModSWfour}

\newcommand{\SpinC}{{\mathrm{Spin}}^c}

\newcommand\Wedge{\Lambda}

\newcommand\Hom{\mathrm{Hom}}

\newcommand\abuts\Rightarrow
\newcommand\Sym{\mathrm{Sym}}

\newcommand\mCP{{\overline{\mathbb{CP}}}^2}

\newcommand\HFpRed{\HFp_{\red}}

\newcommand\relspinc{\underline{\spinc}}

\newcommand\x{\mathbf x}

\newcommand\y{\mathbf y}

\newcommand\ModSphere{\ModFlow\left({\mathbb S}\longrightarrow 
\Sym^{g-1}(\Sigma_{1})\times \Sym^2(\Sigma_{2})\right)}
\newcommand\ModSpheres\ModSphere

\newcommand\CFa{\widehat{CF}}

\newcommand{\red}{\mathrm{red}}

\newcommand\HFp{\HFb}

\newcommand\HFa{\widehat{HF}}
\newcommand\HFb{HF^+}

\newcommand\Mas{\mu}
\newcommand\UnparModSp{\widehat \ModSp}
\newcommand\UnparModFlow\UnparModSp
\newcommand\Mod\ModSp

\newcommand{\cald}{{\mathcal D}}

\newcommand{\spinc}{\mathfrak s}

\newcommand\ModMaps{\mathcal M}
\newcommand\ModSp\ModMaps

\newcommand\alphas{\mbox{\boldmath$\alpha$}}

\newcommand\etas{\mbox{\boldmath$\eta$}}
\newcommand\betas{\mbox{\boldmath$\beta$}}
\newcommand\gammas{\mbox{\boldmath$\gamma$}}
\newcommand\deltas{\mbox{\boldmath$\delta$}}

\newcommand\spincrel\relspinc

\newcommand\MCone{M}

\newcommand\Dual{\mathcal D}
\newcommand\Duality\Dual

\newcommand\Thetaa{\widehat\Theta}
\newcommand\Ga{\widehat G}
\newcommand\fa{\widehat f}

\newcommand\Black{{\mathcal B}}
\newcommand\ModRect{\Mod(\Square)}
\def\Square{\mathchoice\sqr64\sqr64\sqr{2.1}3\sqr{1.5}3}
\def\sqr#1#2{{\vcenter{\vbox{\hrule height.#2pt
	\hbox{\vrule width.#2pt height#1pt \kern#1pt
	\vrule width.#2pt}
	\hrule height.#2pt}}}}
\newcommand\MulFram{I}
\newcommand\DC{\Sigma}
\newcommand\Kh{\mathrm{Kh}}
\newcommand\CKh{\mathrm{CKh}}
\newcommand\CKhred{\widetilde\CKh}
\newcommand\Khred{\widetilde\Kh}
\newcommand\Diagram{\mathcal D}
\newcommand\dkh{\mathfrak d}
\newcommand\Char{\mathrm{Ch}}
\newcommand\Formal{\widehat{\mathbb H}}

\newcommand\Vred{\widetilde V}

\commentable{

\title[{On the Heegaard Floer homology of branched double-covers}]
{On the Heegaard Floer homology of branched double-covers}

\author[Peter Ozsv{\'a}th]{Peter Ozsv\'ath}
\address{Department of
Mathematics, Columbia University, New York 10027 \newline
\indent{\emailfont{petero@math.columbia.edu}}}
\thanks{PSO was supported by NSF grant number DMS 9971950 and a Sloan 
Research Fellowship}

\author[Zolt{\'a}n Szab{\'o}]{Zolt{\'a}n Szab{\'o}} 
\address{Department of
Mathematics, Princeton University, New Jersey 08540 \newline
\indent{\emailfont{szabo@math.princeton.edu}}}}
\thanks{ZSz was supported by NSF grant number DMS 0107792
and a Packard Fellowship}

\begin{document}

\begin{abstract}  
        Let $L\subset S^3$ be a link.  We study the Heegaard Floer
        homology of the branched double-cover $\Sigma(L)$ of $S^3$,
        branched along $L$. When $L$ is an alternating link, $\HFa$ of
        its branched double-cover has a particularly simple form,
        determined entirely by the determinant of the link. For the
        general case, we derive a spectral sequence whose $E^2$ term
        is a suitable variant of Khovanov's homology for the link $L$,
        converging to the Heegaard Floer homology of $\Sigma(L)$.
\end{abstract}


\maketitle

\section{Introduction}

Given a link $L\subset S^3$, we can form its branched double cover, a
new three-manifold which we denote by $\DC(L)$. In this paper, we
study the Heegaard Floer homology of this three-manifold
$\HFa(\DC(L))$ (c.f.~\cite{HolDisk}).

\begin{figure}
\mbox{\vbox{\epsfbox{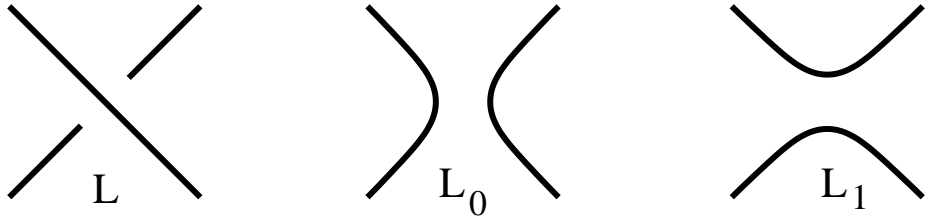}}}
\caption{\label{fig:Skein}
{\bf{Skein moves.}}
Given a link with a crossing as labelled in $L$ above, we have two
``resolutions'' $L_0$ and $L_1$, obtained by replacing the crossing by
the two simplifications pictured above.}
\end{figure}

The starting point for these investigations is a skein exact sequence
which this link invariant $L\mapsto \HFa(\DC(L))$
satisfies. Specifically, fix a projection of $L$, and let $L_0$ and
$L_1$ denote the two resolutions of $L$ at a crossing for the
projection, as illustrated in Figure~\ref{fig:Skein}.  It is a quick
consequence of the surgery long exact sequence for $\HFa$ that for any
link $L\subset S^3$, the groups $\HFa(L)$, $\HFa(L_0)$, and
$\HFa(L_1)$ fit into a long exact sequence
\begin{equation}
\label{eq:DCSkeinExactSeq}
\begin{CD} ...@>>> \HFa(\Sigma(L_0)) @>>> \HFa(\Sigma(L_1)) @>>> \HFa(\Sigma(L))@>>> ...
\end{CD}
\end{equation}

This skein exact sequence leads readily to a complete calculation of
$\HFa(\Sigma(L))$, where $L$ is any alternating link, explained in
Section~\ref{sec:AltLink}. In particular, it is shown there that if
$L$ is a link which admits a connected, alternating projection, then
the rank of $\HFa(\Sigma(L))$ agrees with the number of elements in
$H^2(\Sigma(L);\Z)$, i.e. that $\Sigma(L)$ is what might be called an
``ungraded Heegaard Floer homology lens space'' or, in the terminology
of~\cite{NoteLens}, an $L$-space. (Note that this argument applies to
a larger class of links which includes all alternating links as a
special case.)

When $Y$ is an arbitrary three-manifold, $\HFa(Y)$ has the structure
of a $\Zmod{2}$-graded Abelian group, and that is the structure we
will be concerned with throughout most of this paper.  But in general,
$\HFa(Y)$ also comes with a natural splitting into summands indexed by
$\SpinC$ structures on $Y$~\cite{HolDisk}. Indeed, when $Y$ is a
rational homology three-sphere, the groups are further endowed with an
absolute $\Q$-grading~\cite{AbsGraded}.

By further elaborating on the calculations for $\Sigma(L)$ when $L$ is
alternating, we are able to determine this extra structure explicitly
from the alternating diagram for $L$, as explained in
Section~\ref{sec:AltLink}. As explained in~\cite{AbsGraded} (compare
also~\cite{Froyshov}), this structure gives constraints on the
intersection forms of negative-definite four-manifolds which bound
$\Sigma(L)$. 

Turning back to the case of a general link $L$, it is suggestive to
compare the Exact Sequence~\eqref{eq:DCSkeinExactSeq} with the work of
Khovanov, c.f.~\cite{Khovanov} (for the reader's convenience, we
briefly review the construction in Section~\ref{sec:Khovanov}).
Specifically, Khovanov introduces an invariant for links in $S^3$
whose Euler characteristic, in a suitable sense, is the Jones
polynomial (c.f.~\cite{Jones}, see also~\cite{Kauffman}).  By
construction, his invariants satisfy a ``skein exact sequence''
inspired by the skein relation for the Jones polynomial.  In
particular, just like $\HFa(\Sigma(L))$, Khovanov's invariants fit
into a long exact sequence relating the invariant for a link and its
two resolutions:
\begin{equation}
\label{eq:KhSkeinExactSeq}
\begin{CD} ...@>>> \Kh(r(L_0)) @>>> \Kh(r(L_1)) @>>> \Kh(r(L))@>>> ...,
\end{CD}
\end{equation}
where here $r(L)$ denotes the mirror of $L$. (Note that our
conventions on $L_0$ and $L_1$ are opposite to Khovanov's; this is
why we use the mirror.)
But unlike $\HFa(\Sigma(L))$, Khovanov's theory comes with extra
gradings (which the maps in the exact sequence respect), which allow
one to extract the Jones polynomial from the Betti numbers. (Note that
we are using Khovanov's homology, rather than cohomology, and hence
our conventions on $L_0$ and $L_1$ are opposite to those
of~\cite{Khovanov}.)

The connection between the two link invariants is provided by the
following result.

\begin{theorem}
\label{intro:DCovSpecSeq}
Let $L\subset S^3$ be a link.  There is a spectral sequence whose
$E^2$ term consists of Khovanov's reduced homology of the mirror
of $L$ with
coefficients in $\Zmod{2}$, and which converges to $\HFa(\Sigma(L);\Zmod{2})$.
\end{theorem}

See Section~\ref{sec:DCovSpecSeq} for a precise statement, and also
the proof.  Note that in the above statement, we use here a
``reduced'' version of Khovanov's homology, which he introduced
in~\cite{KhovTwo}, with coefficents in $\Zmod{2}$.

We have the following quick corollary (whose proof is spelled out in
Section~\ref{sec:DCovSpecSeq}):

\begin{cor}
\label{cor:DetKhovBounds}
Let $L\subset S^3$ be a link, and let $\rk \Khred(L)$ denote the rank of its
reduced Khovanov homology with $\Zmod{2}$ coefficients. Then, we have 
the inequality
$$\det(L)\leq \rk_{\Zmod{2}}\HFa(\Sigma(L);\Zmod{2}) \leq \rk_{\Zmod{2}}\Khred(L),$$
where here $\det(L)$ denotes the determinant of the knot (i.e.
th
\end{cor}

Theorem~\ref{intro:DCovSpecSeq} is seen as a consequence of a ``link
surgeries spectral sequence'' established in
Section~\ref{sec:LinkSurgeries}, which holds in a more general setting
(c.f. Theorem~\ref{thm:LinkSpecSeq}). To place this result in context,
recall that if $K\subset Y$ is a framed knot in a three-manifold,
in~\cite{HolDiskTwo}, it is shown that if $Y_0$, $Y_1$ denote the
result of surgeries on $Y$ along $K$ (here, as usual, $Y_0$ denotes surgery
along $K$ in $Y$ with respect to the given framing, while $Y_1$
denotes surgery along $K$ in $Y$ with respect to the framing obtained
by adding a meridian to the given framing), then there is a long exact
sequence relating $\HFa(Y)$, $\HFa(Y_0)$, and $\HFa(Y_1)$, compare
also~\cite{FloerTriangles}. When the knot is replaced by a
multi-component link, the corresponding object is a spectral sequence
relating the various surgeries on the various components of the
link. This spectral sequence, in turn, is established with the help of
the associativity properties of the pseudo-holomorphic polygon
construction, see also~\cite{FOOO},
\cite{QAssoc}, \cite{Silva}.

To establish Theorem~\ref{intro:DCovSpecSeq} we specialize the link
surgeries spectral sequence of Section~\ref{sec:LinkSurgeries} to the
case arising from the branched double cover of a link projection.
Given a projection of $L$, $\Sigma(L)$ comes equipped with a link,
whose components correspond to crossings in the projection, framed so
that surgeries on these components give branched double-covers of the
resolutions of $L$ (this is the topological input for establishing
Equation~\eqref{eq:DCSkeinExactSeq}).  With this said, the key
observation leading to Theorem~\ref{intro:DCovSpecSeq} is the
following.  Consider the the branched double cover of a collection of
unlinks in the plane, connected by cobordisms induced from the
connected sums among (and within) the circles. Applying $\HFa$ (with
coefficients in $\Zmod{2}$) to these objects and morphisms (as
required in the $E^1$ term coming from the link surgery spectral
sequence), one recaptures the $(1+1)$-dimensional topological quantum
field theory which underpins Khovanov's invariants. Armed with this
observation, Theorem~\ref{intro:DCovSpecSeq} follows quickly.

\subsection{Further remarks and speculation}

It is interesting to note that the results from
Section~\ref{sec:AltLink} on non-split alternating links, can be interpreted as
saying that the spectral sequence of Theorem~\ref{intro:DCovSpecSeq}
(c.f. Theorem~\ref{thm:DCovSpecSeq} for the precise statement)
collapses at the $E^2$ stage. (Note that it
is clear from the more precise statement that if $L$ is an
$n$-crossing link, then the spectral sequence always collapses after
the $E^n$ stage.)

A rather striking example where there
are non-trivial differentials beyond the $E^2$ stage
is illustrated for torus knots. For example,
let $T_{p,q}$ denotes the $(p,q)$ torus knot. When both
$p$ and $q$ are odd, $\Sigma(T_{p,q})$ is the Brieskorn homology sphere with
multiplicities $2$, $p$, and $q$. In particular,
$\Sigma(T_{3,5})$ is the Poincar\'e homology sphere, so
$\HFa(\Sigma(T_{3,5};\Zmod{2}))\cong \Zmod{2}$ (c.f.~\cite{AbsGraded}
or~\cite{SomePlumb}), while its reduced
Khovanov homology clearly has larger rank, as its Jones polynomial 
has three non-zero coefficients.

Results from this paper raise a number of further questions, which
further link Khovanov's essentially combinatorial theory with problems
involving holomorphic disks.

As a first point, observe that the link surgeries spectral sequence
works with $\Z$, rather than $\Zmod{2}$ coefficients. This suggests an
$E^2$ term whose $\Zmod{2}$ reduction agrees with Khovanov's reduced
theory, but which differs from the sign conventions as defined by
Khovanov. It would be interesting to pin down the new sign conventions
of this theory, not only from the point of view of applications to
Heegaard Floer homology (i.e. to give information about $\HFa$ over
$\Z$ of the branched cover), but also from the point of view of
Khovanov's theory, as it would give a link invariant with $\Z$-coefficients
whose Euler
characteristic is the normalized Jones polynomial.

In another direction, it is reasonable to expect that the induced
filtered quasi-isomorphism type associated to the branched double
cover spectral sequence from Theorem~\ref{intro:DCovSpecSeq} is also a
link invariant, i.e. that is independent of the projection used in its
definition. This would give a countable sequence of link invariants,
starting with Khovanov's homology, and ending with $\HFa$ of the
branched double cover.

This also raises the question of finding a combinatorial
description of the higher differentials for the spectral
sequence. Although finding a combinatorial description of the Heegaard
Floer homology in general is a very interesting, if difficult problem,
it is perhaps easier when one specializes to the case of
branched double covers of links in the three-sphere.

Another question concerns naturality properties of Khovanov's
homology. On the one hand, Khovanov~\cite{Khovanov} shows that a knot
cobordism $X$ from $L_1$ to $L_2$ induces a (combinatorially defined)
map between Khovanov homologies. Now, the branched double-cover of
$X$ inside $[1,2]\times S^3$
is a four-manifold $\Sigma(X)$ which gives a cobordism from
$\Sigma(L_1)$ to $\Sigma(L_2)$, and correspondingly induces a map on
$\HFa$ (c.f.~\cite{HolDiskFour}), defined by counting holomorphic
triangles. This map, in general, can be quite difficult to compute. It
is reasonable to expect that there is a well-defined map between the
filtered complexes which give rise to 
Theorem~\ref{intro:DCovSpecSeq} and hence between spectral sequences
which at the $E^2$ stage induces Khovanov's map,
and  at the $E^\infty$ stage induces the map on $\HFa$ induced by
$\Sigma(X)$.

\subsection{Organization}
The skein exact sequence for $\HFa(\Sigma(L))$ is established in
Section~\ref{sec:BranchedCovers}; the results for alternating links
(with a sample calculation) are explained in
Section~\ref{sec:AltLink}.  The link surgeries spectral sequence is
established in Section~\ref{sec:LinkSurgeries} (note that this general
result applies not only to branched double covers, considered in the
rest of the paper). In Section~\ref{sec:Khovanov}, we review
Khovanov's link invariant (with $\Zmod{2}$ coefficients), setting up
the notation for Section~\ref{sec:DCovSpecSeq}, where we establish
the precise form of Theorem~\ref{intro:DCovSpecSeq}.

\subsection{Acknowledgements}
Similarities between Khovanov's reduced theory and $\HFa$ of a
branched double cover were first noticed by Jacob Rasmussen and
Ciprian Manolescu.  We would like to thank Jacob Rasmussen for some
very helpful discussions on Khovanov's homology, and also Paul Seidel
for some very helpful discussions on homological algebra.

\section{Skein moves and branched double covers}  
\label{sec:BranchedCovers}

Let $K$ be a framed knot in a three-manifold $Y$ (i.e. a knot with a
choice of longitude $\lambda$). Let $Y_0=Y_0(K)$ denote the
three-manifold obtained from $\lambda$-framed surgery on $Y$ along
$K$, and let $Y_1=Y_1(K)$ denote the three-manifold obtained from
$\lambda+\mu$-framed surgery on $Y$ along $K$ (where here $\mu$
denotes the canonical meridian for the knot $K$). We call the ordered 
triple $(Y,Y_0,Y_1)$ a {\em triad} of three-manifolds.

This relationship between $Y$, $Y_0$, and $Y_1$ is symmetric under a
cyclic permutation of the the three three-manifolds. Indeed, it is not
difficult to see that $(Y,Y_0,Y_1)$ fit into a triad if and only if
there is a single oriented three-manifold $M$ with torus boundary, and
three simple, closed curves $\gamma$, $\gamma_0$, and $\gamma_1$ in
$\partial M$ with 
\begin{equation}
\label{eq:TriadRelation}
\#(\gamma\cap \gamma_0)=\#(\gamma_0 \cap \gamma_1)
=\#(\gamma_1\cap \gamma) = -1
\end{equation} (where here the algebraic intersection
number is calculated in $\partial M$, oriented as the boundary of
$M$), so that $Y$ resp. $Y_0$ resp. $Y_1$ are obtained from $M$ by
attaching a solid torus along the boundary with meridian $\gamma$
resp. $\gamma_0$ resp. $\gamma_1$.

In~\cite{HolDiskTwo}, we established a long exact sequence connecting
$\HFa$ for any three three-manifolds which fit into a triad: $$
\begin{CD}
...@>>>\HFa(Y)@>>>\HFa(Y_0)@>>>\HFa(Y_1)@>>>...
\end{CD}
$$

The skein exact sequence for $\HFa(\Sigma(L))$
(Equation~\eqref{eq:DCSkeinExactSeq}) follows readily:

\begin{prop}
\label{prop:DCSkeinExactSeq}
Fix a crossing for a projection of a link $L\subset S^3$,
and $L_0$ and $L_1$ be the two resolutions of that crossing as in Figure~\ref{fig:Skein}. Then, the three-manifolds $(\Sigma(L), \Sigma(L_0), \Sigma(L_1))$
form a triad. In particular, there is an induced long exact sequence
$$\begin{CD}
...@>>> \HFa(\DC(L_0)) @>>> \HFa(\DC(L_1)) @>>> \HFa(\DC(L)) @>>> ...
\end{CD}$$
\end{prop}

\begin{proof}
Fix a sphere $S$ meeting the link $L$ in four points, containing a
ball $B$ which contains two arcs of $L$, and in whose complement $L$,
$L_0$, and $L_1$ agree. Clearly, letting $M$ be the branched
double-cover of $S^3-B$ branched along $L-(L\cap B)$, we see that
$\Sigma(L)$, $\Sigma(L_0)$, and $\Sigma(L_1)$ are all obtained from
$M$ by attaching the branched double cover of $B$ branched along two
arcs. 

Now, it is easy to see that the branched double-cover of $B$ branched
along two standard, unknotted arcs is a solid torus. Indeed, a
meridian for this solid torus can be realized by pushing either of the
two arcs out to the boundary, and taking its branched
double-cover. Thus, letting letting
$\gamma$, $\gamma_0$, and $\gamma_1$
denote curves obtained by pushing arcs out into the boundary torus
for the $L$ and its resolutions $L_0$ and $L_1$, it is straightforward to
verify that these curves satisfy Equation~\eqref{eq:TriadRelation}.

Thus, $(\Sigma(L), \Sigma(L_0),\Sigma(L_1))$ form a triad of
three-manifolds.  The exact sequence now is a direct consequence of
the aformentioned surgery long exact sequence
(Theorem~\ref{HolDiskTwo:thm:GeneralSurgery} of~\cite{HolDiskTwo};
see Theorem~\ref{thm:KnotCase} below for another proof).
\end{proof}

\begin{figure}
\mbox{\vbox{\epsfbox{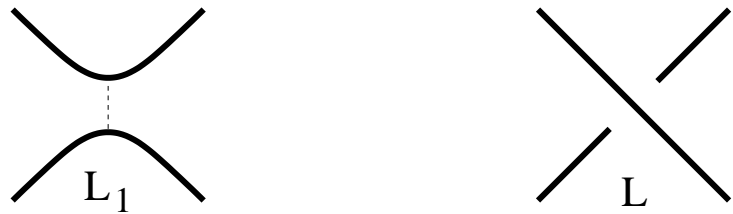}}}
\caption{\label{fig:CoveringTwoHandle}
{\bf{Obtaining $\Sigma(L)$ from $\Sigma(L_1)$.}} 
The three-manifold $\Sigma(L)$ (the corresponding
to the branched double cover of a link with
the a crossing as illustrated on the right)
is obtained from $\Sigma(L_1)$
by surgery on the knot obtained as a branched double cover 
of the dashed arc indicated in the picture on the left.}
\end{figure}

In particular, we have seen that $\Sigma(L)$ is obtained as surgery on
a knot in $\Sigma(L_1)$. This knot can be explicitly seen as the
branched double cover of a standard arc inside the three-ball $B$
containing the two resolved arcs in $L_1$. In turn, this arc can be
pictured in a knot projection of $L_1$ as an arc $A$ which meets $L_1$
in exactly two points, both of which are on the boundary of $A$, and
which connect the two resolved strands in $L_1$, as pictured in
Figure~\ref{fig:CoveringTwoHandle}.

\section{Alternating links}
\label{sec:AltLink}

Let $Y$ be an oriented three-manifold. Let $|H^2(Y;\Z)|$ denote
the number of elements in $H^2(Y;\Z)$ provided that $b_1(Y)=0$,
and let $|H^2(Y;\Z)|=0$ if $b_1(Y)>0$.
Now, if $L$ is a link in $S^3$,  the determinant of $L$ is defined
by $\det(L)=|\Delta_{L}(-1)|$, where here $\Delta_L(T)$ denotes the
Alexander polynomial of $L$. It is well-known (see for
example~\cite{Lickorish}) that $\det(L)=|H^2(\Sigma(L);\Z)|$.

Recall that the Euler characteristic of $\HFa(Y)$ is given by
$|H^2(Y;\Z)|$ (c.f. Proposition~\ref{HolDiskTwo:prop:EulerHFa} of~\cite{HolDiskTwo}); in particular,
$|H^2(Y;\Z)|\leq \rk \HFa(Y)$.  Three-manifolds with $b_1(Y)=0$ for
which $|H^2(Y;\Z)|=\rk \HFa(Y)$ are called $L$-spaces
(c.f.~\cite{NoteLens}).  This special class of three-manifolds is
closed under connected sums and includes all lens spaces and, more
generally, all Seifert fibered spaces with finite fundamental group;
other examples are given in~\cite{SomePlumb} and~\cite{NoteLens}.  We
will prove that if $L$ is a non-split, alternating link, then
$\Sigma(L)$ is an $L$-space. Indeed, the class of links we work with
here is wider. To this end, we have the following:

\begin{defn}
The set ${\mathcal Q}$ of {\em quasi-alternating links} is
the smallest set of links which satisfies the following properties:
\begin{enumerate}
\item the unknot is in ${\mathcal Q}$
\item the set ${\mathcal Q}$ is closed under the following operation.
Suppose $L$ is any link which admits a projection
with a crossing with the following properties:
\begin{itemize}
\item both resolutions $L_0, L_1\in {\mathcal Q}$, 
\item $\det(L_0), \det(L_1)\neq 0$, 
\item $\det(L)=\det(L_0)+\det(L_1)$;
\end{itemize}
then $L\in {\mathcal Q}$.
\end{enumerate}
\end{defn}

Note that quasi-alternating in this sense is different from the notion
of almost-alternating, which appears in the literature (c.f.~\cite{AlmostAlternating}).

\begin{lemma}
Every link which admits a connected, alternating
projection  is quasi-alternating.
\end{lemma}

\begin{proof}
Recall that a complement of a knot projection in the plane admits a
checkerboard coloring.  The collection of black regions can be given
the structure of a planar graph $\Black(L)$, whose vertices correspond to
black regions and edges correspond to vertices which are corners of
pairs of black regions.  It is a classical result~\cite{Crowell}
that if $L$ admits an alternating projection, then the
determinant of $L$ is the total number of maximal subtrees of the
black graph of $L$. To fix orientation conventions, when coloring an
alternating link, we always use the coloring scheme indicated in
Figure~\ref{fig:ColorConventions}. (Although this choice is irrelevant
to the present discussion, it becomes vital in pinning down signs in 
our later statements, especially in Theorem~\ref{thm:AbsGradeAlt}
and Equation~\eqref{eq:GoLi} below.)

\begin{figure}
\mbox{\vbox{\epsfbox{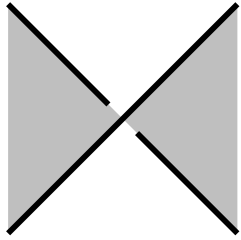}}}
\caption{\label{fig:ColorConventions}
{\bf{Coloring conventions for alternating knots.}} 
We adopt the pictured convention when coloring an alternating projection.}
\end{figure}

We now induct on the determinant of the link.  In the basic case where
the determinant is one, it follows at once that there is only one
maximal subtree, and hence that the knot is the unknot.  

For the inductive step,
it is easy to see that for a reduced alternating projection of $L$, if
we choose any crossing $x$, both resolutions $L_0$ and $L_1$ at $x$
are connected, alternating projections of links. Moreover, it is easy
to see that $\det(L)=\det(L_0)+\det(L_1)$: maximal subtrees of the
black graph of $L$ which contain, resp. do not contain, the edge
corresponding to $x$ are in one-to-one correspondence with the maximal
subtrees of the black graph of $L_i$, resp $L_j$, where here $i,j\in
\{0,1\}$ and $i\neq j$. Thus, by the inductive hypothesis, 
the theorem has been established for both $L_0$ and $L_1$; and hence,
the inductive step follows.
\end{proof}

\begin{figure}
\mbox{\vbox{\epsfbox{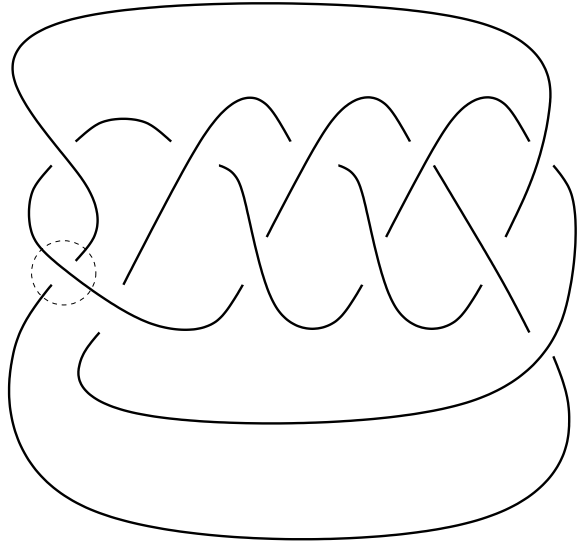}}}
\caption{\label{fig:Nine47}
{\bf A quasi-alternating, but not alternating, knot.}
The pictured knot $9_{47}$ is quasi alternating: its determinant is $29$,
and if we resolve the indicated crossing either way, we obtain
(non-split) alternating links with determinants $5$ and $24$.}
\end{figure}

Of course, there are quasi-alternating links which are not alternating.
For a picture of one, see Figure~\ref{fig:Nine47}.

\begin{prop}
\label{prop:AltLink}
If $L$ is a quasi-alternating link, $\Sigma(L)$ is an $L$-space, i.e.
$$\HFa(\Sigma(L))\cong
\Z^{\det(L)}.$$
\end{prop}

\begin{proof}
The proposition is now established by induction on the determinant of $L$.
In the basic case where the determinant is one, it follows at once that
there is only one maximal subtree, and hence that the knot is the unknot,
so $\HFa(\Sigma(L))=\HFa(S^3)\cong \Z$.

The bound $\det(L)\leq \rk \HFa(\Sigma(L))$ combined with the long
exact sequence Proposition~\ref{prop:DCSkeinExactSeq} readily provides
the inductive step (c.f. Proposition~\ref{NoteLens:prop:LSpaces}
of~\cite{NoteLens}).
\end{proof}

When $L$ is quasi-alternating, it is not difficult to calculate the
absolute gradings of $\HFa(\Sigma(L))$.  To this end, we find it
convenient to restrict attention to the case where $L$ is a 
non-split, alternating link.

Let $L$ be a link with a connected, alternating projection, and choose
a maximal subtree $T$ of the black graph $\Black(L)$, and let
$\{e_i\}_{i=1}^m$ denote the edges in $Z_T=\Black(L)-T$. Let $V$ denote the
lattice generated by these edges. We can equip $V$ with a bilinear
form $$Q\colon V\otimes V \longrightarrow \Z$$ 
as follows.
Choose orientations for each edge $e_i\in Z_T$, let $C_i$ denote
the oriented circuit in $T\cup \{e_i\}$; 
and if $X$ is any subgraph of $\Black(L)$,
let $E(X)$ denote the number of edges in $X$. Note that the
orientation on $e_i$ induces an orientation on the circuit $C_i$.
Given a pair of
distinct edges $e_i, e_j\in Z_T$ with the property that $C_i\cap
C_j\neq \emptyset$, we let
$$Q(e_i\otimes e_j)= \epsilon(i,j) \cm E(C_i\cap C_j),$$
where here $\epsilon(i,j)$ is given by
$$\epsilon(i,j)= 
\left\{\begin{array}{ll}
+1 &{\text{if the orientation on $C_i\cap C_j$ induced from $C_i$}} \\
& {\text{is opposite to the one induced from $C_j$;}} \\
-1 & {\text{otherwise.}}
\end{array}
\right.
$$
In particular, $Q(e_i\otimes e_i)=-E(C_i)$.

We will realize $Q$ as the intersection
form of a four-manifold which bounds $\Sigma(L)$ in the proof of
Theorem~\ref{thm:AbsGradeAlt} below.

A characteristic vector for $M$ is a vector in the lattice $K\in
V^*$ with $\langle K, v\rangle + Q(v,v)\equiv 0 \pmod{2}$ for
each $v\in V$. Two characteristic vectors $K$ and $K'$ are said to 
be equivalent if $K-K'=2Q(v\otimes \cdot)$ for some $v\in V$.

\begin{theorem}
\label{thm:AbsGradeAlt}
There is an identification $i$ 
equivalence classes of characteristic vectors for $Q$
with $\SpinC$ structures over $\Sigma(L)$.
Moreover, given an equivalence class of characteristic
vectors ${\Xi}$, 
$\HFa(\Sigma(L),i({\Xi}))\cong \Z$ is supported in dimension
$$d(\Xi)=\frac{\max_{K\in{\Xi}}K^2+b}{4},$$
where here $K^2$ is the length of $K$ with respect to the inner product on
$V^*$ induced from $Q$, and $b$ is the number of edges in $Z_T$
(or, more invariantly, the rank of $H_1(\Black(L);\Z)$).
\end{theorem}

\begin{remark}
We emphasize that we are using the coloring conventions pictured in
Figure~\ref{fig:ColorConventions}, which breaks the apparent symmetry
between the ``white'' and ``black'' graphs. In fact, using the white
graph in place of the black graph to construct the form analogous to
$Q$, it is not difficult to see that we obtain the dimensions of the
generators for $\HFa(-\Sigma(L))$, whose sign is opposite to those for
$\HFa(\Sigma(L))$.
\end{remark}

We break the proof into several pieces. First, we construct a
four-manifold $X_L$ which bounds $\Sigma(L)$.  To construct $X_L$, fix
a projection for the link $L$, and let $n$ denote its number of
crossings. If we form $1$-resolutions at each intersection, we obtain
a $k$-component unlink. The branched double cover of this manifold is
$Y_0=\#^{k-1}(S^2\times S^1)$.  Attaching one two-handle for each
crossing to ``unresolve'' the crossing (as in
Proposition~\ref{prop:DCSkeinExactSeq}), we obtain a cobordism from
$Y_0$ to $\Sigma(L)$. Indeed, by filling $Y_0$ by the boundary
connected sum of $k-1$ copies of $B^3\times S^1$, we obtain a
four-manifold $X_L$ which bounds $Y$. 

\begin{lemma}
\label{lemma:IntFormFour}
If $L$ is a non-split alternating link, the four-manifold $X_L$
described above is negative-definite. Indeed, there is an
identification of the form $Q$ on the vector space $V$ described above
with the intersection form on the two-dimensional homology of $X_L$.
\end{lemma}

\begin{proof}
As we have described it, $X_L$ is built from one zero-handle, $k-1$
one-handles, and $n$ two-handles. In fact, the tree $T$ specifies
$k-1$ two-handles which cancel the one-handles; i.e. after attaching the
two-handles from the tree, we obtain the branched double cover of a
single unknot, which is $S^3$. Now, $X_L$ is obtained from the four-ball
by surgery on a link in $S^3$ (the branched double cover of the unknot)
whose components correspond to the remaining edges 
in $\Black(L)-T$ (i.e. the link components
are the branched double covers of the arcs with boundary 
in the unknot, which are associated
to the edges in $\Black(L)-T$).
We claim that a choice of orientation on each edge
$e_i$ simultaneously orients all the components of this link, up to 
an overall sign.

To see this, we proceed as follows. Let $U$ denote the unknot as
specified by the tree $T$. Let $\phi$ be a vector field normal to $U$
which is orthogonal to the kernel of the projection map used in
describing the knot projection. This vector field $\phi$, of course,
specifies the blackboard framing of $U$.  The vector field $\phi$ has
two possible lifts in the branched double cover of the unknot
(in the sense that there are two lifts in the branched double cover
of the knot obtained by displacing
$U$ by $\phi$).  Choose
one, and denote it ${\widetilde \phi}$ (while the other is denoted
${\widetilde\phi}'$). Then, the knot corresponding to $e_i$ (thought
of as an arc connecting $x$ to $y$ in the unknot) is oriented so that
its tangent vector at $x$ agrees with ${\widetilde \phi}_x$ (as
opposed to ${\widetilde \phi}'_x$). We denote the oriented knot
associated to $e_i$ with its orientation by $k_i$.  (Note that the
other lift of the blackboard framing has the effect of reversing the
induced orinetations on all the knots $k_i$ simulatneously.)

Next, we argue that the intersection form of $X_L$ is
negative-definite. We prove this by induction on the number of
crossings.  The basic case is obvious. Next, recall that
$$\Big|H^2(\Sigma(L);\Z)\Big|=\Big|H^2(\Sigma(L_0);\Z)\Big|+
\Big|H^2(\Sigma(L_1);\Z)\Big|,$$  
so it readily follows that the two-handle from $\Sigma(L_1)$ to
$\Sigma(L)$ (and also the one from $\Sigma(L)$ to $\Sigma(L_0)$) is
negative-definite.  Now, it is easily seen that $X_L$ is obtained by
attaching this negative-definite two-handle to $X_{L_1}$. 

We show that the intersection form on $H_2(X_L;\Z)$ is given by
$Q$. To this end, observe that if we attach $m$ of the remaining two
handles $\{e_i\}_{i=1}^m$, to $S^3$, the number of elements in $H^2$
of the boundary three-manifold (with $\Z$ coefficients) is given by
the determinant of the matrix $(Q(e_i\otimes e_j))_{i,j\in
\{1,...,n\}}$, which in turn is obtained from the number of maximal
subtrees in $T\cup\{e_i\}_{i=1,...,m}$.

In particular, since $X_L$ is negative-definite, it follows at once
that if we choose a basis for $H_2(X_L;\Z)$ given by the two-handles
in $X_L$ (with any set of orientations), then if $[e_i]$ the homology
class corresponding to an edge $e_i$, then
$\#[e_i]\cap[e_i]=-E(C_{i})$ (since the number of maximal subtrees of
a circuit is the length of the circuit, and the sign is forced by the
negative-definiteness). 

\begin{figure}
\mbox{\vbox{\epsfbox{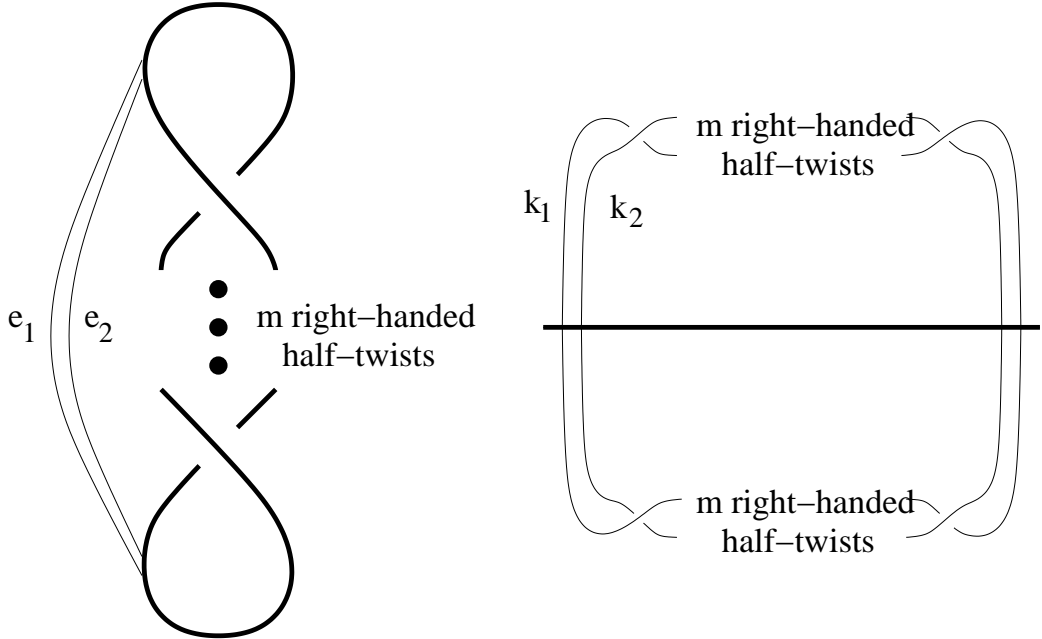}}}
\caption{\label{fig:StandardPicture}
{\bf{Standard picture for two crossings.}}  If $e_1$ and $e_2$
correspond to two closed circuits with $m=E(C_1\cap C_2)$ edges in
common, then the unknot corresponding to the tree $T$, together with
the two arcs associated to $e_1$ and $e_2$, is isotopic to the picture
on the left (which depicts a projection of the unknot with $m$
right-handed half-twists in it, of which two have already been drawn,
together with two arcs which meet the
unknot in the specified manner). Passing to the branched double cover
of the unknot (which in turn is best visualized by unwinding the
twists on the unknot, at the expense of twisting the arcs
corresponding to $e_1$ and $e_2$, and then taking the branched double
cover of these edges to obtain knots $k_1$ and $k_2$, we obtain the picture shown on the right, where the solid line
indicates the branched locus.}
\end{figure}

When $i\neq j$, $\#[e_i]\cap [e_j]$ is given by the linking number of
$k_i$ with $k_j$. This in turn is calculated in a model case:
consider the unknot corresponding to the tree $T$, together with the
two arcs corresponding $e_i$ and $e_j$. This is easily seen to be
isotopic to an unknot with two arcs attached, in a manner which has a
standard projection depending only on the integer $E(C_i\cap C_j)$, as
pictured in Figure~\ref{fig:StandardPicture}.  Unwinding the unknot
and taking the branched double cover, we see that the branched double
covers of the original arcs become circles which are linked
$|E(C_i\cap C_j)|$ times. It is then straightforward to see that the
sign of this linking number is the one stated (once we choose a lift
of the blackboard framing for the unknot). 
\end{proof}

\vskip.3cm
\noindent{\bf{Proof of Theorem~\ref{thm:AbsGradeAlt}.}}
With Lemma~\ref{lemma:IntFormFour}, the proof of the theorem now
follows along the lines of Section~\ref{SomePlumb:sec:Proof}
of~\cite{SomePlumb}. We sketch here the main points.
Let $\Char(X_L)$ denote the set of
characteristic vectors for the intersection form $H^2(X_L;\Z)$.  We
write $K\sim K'$ if there is an element $v\in H^2(X_L,\partial X_L)$
with the property that $K=K'+2v$.  Next (compare~\cite{SomePlumb}),
consider the subgroup $$\Formal(X_L) \subset \Hom(\Char(X_L),\Z),$$
consisting of maps $\phi$ with the properties that
\begin{itemize}
\item $\phi(K)=\phi(K')$ if $K\sim K'$ and $Q(K,K)=Q(K',K')$
\item $\phi(K)=0$ if there is some $K'\sim K$ with $Q(K',K')> Q(K,K)$.
\end{itemize}
Viewing $X_L$ as a cobordism from $-\Sigma(L)$ to $S^3$, we obtain a 
naturally induced map
(c.f.~\cite{HolDiskFour})
$$T_{X_L}\colon \HFa(-\Sigma(L)) \longrightarrow \Hom(\Char(X_L),\Z),$$
in view of the fact that $\HFa(S^3)\cong \Z$. 

Unless the diagram of $L$ represents the unknot, we can always find
a double-point $p$ whose two resolutions are connected diagrams. 
This gives the following commutative diagram:
\begin{equation}
\label{eq:FormalDiag}
\begin{CD}
0@>>> \HFa(-\Sigma(L_0)) @>>>\HFa(-\Sigma(L))
@>>>\HFa(-\Sigma(L_1))@>>> 0 \\ 
&& @V{T_{X_{L_0}\#\mCP}}VV  @V{T_L}VV  @V{T_{L_1}}VV \\ 
&&\Hom(\Char(X_{L_0}\#\mCP),\Z) @>{\mathbb A}>>\Hom(\Char(X_L),\Z)
@>{\mathbb B}>>\Hom(\Char(X_{L_1}),\Z)), \\
\end{CD}
\end{equation}
where the top row is exact, the squares commute, the maps ${\mathbb A}$
and ${\mathbb B}$ are given by
\begin{eqnarray*}
{\mathbb A}(\phi_0)(K)&=& \sum_{\{K_0\in\Char(X_{L_0}\#\mCP)\big|K_0|_{H^2(X_L;\Z)}=K\}} \phi_0(K_0) \\
{\mathbb B}(\phi)(K)&=& \sum_{\{K_1\in\Char(X_{L_1})\big|K_1|_{H^2(X_{L_1};\Z)}=K\}} \phi(K_1). \\
\end{eqnarray*}
A straightforward induction on the number of crossings in the diagram
shows that the image of $T_{X_L}$ is contained in $\Formal(X_L)$.  The
sphere with square $-1$ contained in the composite cobordism from
$\Sigma(L_0)$ to $\Sigma(L_1)$ through $\Sigma(L)$ is used to show
that ${\mathbb B}\circ {\mathbb A}=0$, and also that ${\mathbb A}$ is
injective. Straightforward homological algebra then shows that $T_L$
is an isomorphism, again, by induction on the number of crossings,
together with Diagram~\ref{eq:FormalDiag}, and an identification
$\Formal(X_{L_0}) \cong \Formal(X_{L_0}\#\mCP)$.
(For a more detailed argument establishing an analogous result, 
see the proof of Lemma~\ref{SomePlumb:lemma:NearlyExact} of~\cite{SomePlumb}.)

Endow $\Formal(X_L)$ with a grading, by declaring an element to be
homogeneous of degree $d$ if it supported on those $K\in\Char(X_L)$
with $$-\left(\frac{K^2+\rk H^2(X_L)}{4}\right)=d.$$ Clearly, $T_{X_L}$ carries
$\HFa_d(-\Sigma(L))$ to $\Formal_d(X_L)$ (c.f.~\cite{AbsGraded}).
Since $\HFa_d(-\Sigma(L))\cong \HFa^{-d}(\Sigma(L))$ (c.f.~\cite{AbsGraded}),
the result now follows. \qed
\vskip.3cm

Note that the long exact sequence can be pushed slightly further than
we have done in the above discussion. For example, recall that if $W$
is a cobordism between two spaces $L$-spaces with $b_2^+>0$, then the
induced map on $\HFp$ is trivial (c.f.~\cite{HolDiskFour}). This gives
at once the result that if $L$ differs from an alternating (or indeed
quasi-alternating) knot by at a single crossing, then all the elements
of $\HFpRed(\Sigma(L))$ have the same $\Zmod{2}$-grading.  Indeed, the
map induced by a two-handle from $\Sigma(L_0)$ to $\Sigma(L_1)$, where
$L_0$ and $L_1$ are both quasi-alternating, determined purely by
homological information. This can be used to give information about
the Heegaard Floer homology of $\Sigma(L)$ when its two resolutions
are quasi-alternating. We do not pursue this any further here,
contenting ourselves instead with a sample calculation illustrating
Theorem~\ref{thm:AbsGradeAlt}.

\subsection{An example: $9_{40}$}

\begin{figure}
\mbox{\vbox{\epsfbox{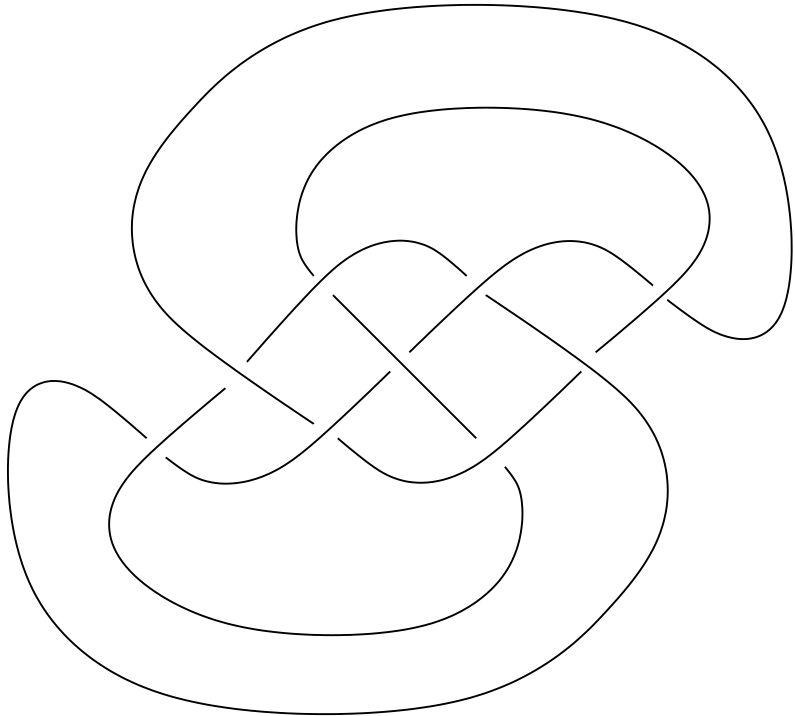}}}
\caption{\label{fig:Nine40}
{\bf The knot $9_{40}$.}.}
\end{figure}

To illustrate Theorem~\ref{thm:AbsGradeAlt}, we calculate
$\HFa(\Sigma(L))$ where $L$ is the alternating knot with nine
crossings $9_{40}$, pictured in Figure~\ref{fig:Nine40}.  

\begin{figure}
\mbox{\vbox{\epsfbox{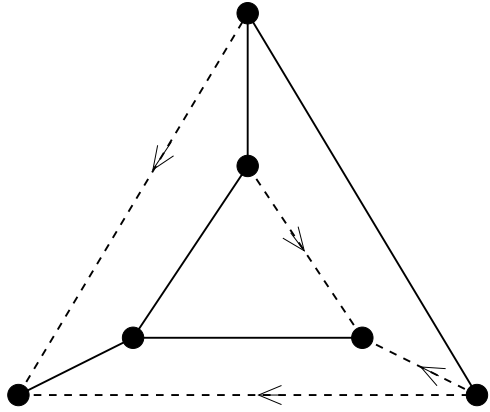}}}
\caption{\label{fig:Nine40w}
{\bf Black graph for the knot $9_{40}$.} 
All edges (including those which are dashed) are inlcuded in the graph.
The solid edges constitute the tree $T$ used for the matrix given in the text.
The dashed edges, when oriented, give rise to the matrix in the text.}
\end{figure}

The black graph of this knot is illustrated in
Figure~\ref{fig:Nine40w}.  Using as our base tree $T$ the solid edges
pictured in the figure, and the orientations of the remaining edges
indicated, the intersection form of $X_L$ takes the form $$Q_L=
\left(
\begin{array}{rrrr}
-3	&-2	&-1	&-1\\
-2	&-5	&-2	&-3\\
-1	&-2	&-4	&-3\\
-1	&-3	&-3	&-5
\end{array}
\right).$$

It is a straightforward if tedious 
matter to find the maximal lengths of the
characteristic vectors for $Q$ in its equivalence classes. Note that
this is a finite search: it is easy to see that all maximal
characteristic vectors have the property that $|\langle
K,v\rangle|\leq |Q(v\otimes v)|$, and hence to determine the absolute
gradings of the generators of $\HFa(\Sigma(L))$
(this and the further calculations in this section were all done
with the help of Mathematica~\cite{Mathematica}). We display the
results below. The numbers are ordered as suggested by the group
structure of $H^2(\Sigma(L);\Z)\cong \Zmod{5}\oplus\Zmod{15}$;
i.e. having chosen such an isomorphism, we have a naturally induced
identification $\SpinC(\Sigma(L))\cong \Zmod{15}\oplus\Zmod{5}$
(where we choose as the origin the spin structure on $\Sigma(L)$; since
$H_1$ has no two-torsion, this structure is uniquely determined); i.e.
the element in the $i^{th}$ row (counting from $0$ to $4$)
and $j^{th}$ column (counting from $0$ to $14$)
is the absolute grading 
of the element in the $\SpinC$ structure corresponding to $(i,j)\in \Zmod{5}\oplus \Zmod{15}$.

$$
\begin{array}{rrrrrrrrrrrrrrrrrrrrrrrrr}
  - \frac{1}{2}&  \frac{13}{30} &
    - \frac{23}{30}&  - \frac{9}{10}&  
    - \frac{11}{30}&  \frac{1}{30}&- \frac{1}{10}&  
    \frac{1}{30}&\frac{1}{30}&- \frac{1}{10}&  
    - \frac{11}{30}&  - \frac{23}{30}&  
    - \frac{9}{10}&  \frac{13}{30}&
    - \frac{11}{30}  \\ \\
   \frac{7}{10}&- \frac{11}{30}&  \frac{13}{30}&
    - \frac{1}{10}&  \frac{13}{30}&\frac{5}{6}&\frac{3}{10}&
    \frac{13}{30}&\frac{13}{30}&- \frac{1}{10}&  
    - \frac{11}{30}&  - \frac{23}{30}&  
    \frac{7}{10}&\frac{1}{30}&- \frac{23}{30}  \\ \\ 
   \frac{3}{10}&- \frac{23}{30}&  \frac{1}{30}&
    - \frac{9}{10}&  - \frac{11}{30}&  \frac{1}{30}&
    - \frac{9}{10}&  - \frac{23}{30}&  
    - \frac{23}{30}&  \frac{3}{10}&\frac{1}{30}&
    - \frac{11}{30}&  \frac{7}{10}&\frac{1}{30}&
    - \frac{23}{30}  \\ \\ 
   \frac{3}{10}&- \frac{23}{30}&  \frac{1}{30}&\frac{7}{10}&
    - \frac{23}{30}&  - \frac{11}{30}&  
    \frac{3}{10}&\frac{13}{30}&\frac{13}{30}&- \frac{9}{10}&  
    \frac{5}{6}&\frac{13}{30}&- \frac{9}{10}&  \frac{13}{30}&
    - \frac{11}{30}  \\ \\ 
   \frac{7}{10}&- \frac{11}{30}&  \frac{13}{30}&\frac{7}{10}&
    - \frac{23}{30}&  - \frac{11}{30}&  
    - \frac{1}{10}&  \frac{1}{30}&\frac{1}{30}&\frac{3}{10}&
    \frac{1}{30}&- \frac{11}{30}&  - \frac{1}{10}&  
    - \frac{23}{30}&  \frac{13}{30}  \\   \\
\end{array}
$$

\section{The link surgeries spectral sequence}
\label{sec:LinkSurgeries}

In this section, we turn our attention away from branched
double-covers, and consider the case of a general three-manifold $Y$.
Our aim here is to describe a generalization of the surgery long exact
sequence for the case of multi-component links in $Y$. In the course
of making this generalization, we give a quick (and slightly stronger)
proof of the long exact sequence based on associativity properties of
the holomorphic polygon construction, combined with some homological
algebra discussed in Subsection~\ref{sec:MCone}.  But first, we
introduce some notation.

Let $L=K_1\cup...\cup K_\ell$ be an $\ell$-component, framed link in a
three-manifold $Y$.  A ``multi-framing'' is a vector
$\MulFram=(m_1,...,m_\ell)$, where each $m_i\in \{0,1,\infty\}$. For a
multi-framing, there is a three-manifold $Y(\MulFram)$, which is
obtained from $Y$ by performing $m_i$-framed surgery on the component
$K_i$ for $i=1,...,n$. As usual, when $m_i=\infty$, this means no
surgery, $m_i=0$ this means $\lambda_i$-framed surgery, and when
$m_i=1$, this is surgery of type $\mu_i+\lambda_i$.

We give the set $\{0,1,\infty\}^\ell$ the lexicographical ordering
(with the understanding that $0<1<\infty$). If
$I\in\{0,1,\infty\}^\ell$, we call $I'$ an {\em immediate successor} of
$I$ if $I=(m_1,...,m_\ell)$ and $I'=(m_1',...,m_\ell')$ if there is
some $j$ so that for all $i\neq j$, $m_i=m_i'$, while $m_j<m_j'$,
excluding the case where $m_j=0$ and $m_j'=\infty$. Clearly, if $I'$ is
an immediate successor of $I$, there is a corresponding map on Floer homology 
$$\Ga_{I<I'}\colon \HFa(Y) \longrightarrow \HFa(Y')$$
associated to the single two-handle addition.

We now state the link surgeries spectral sequence alluded to in the
introduction:

\begin{theorem}
\label{thm:LinkSpecSeq}
There is a spectral sequence whose $E^1$ term is $\bigoplus_{I\in
\{0,1\}^\ell}\HFa(Y(I))$, all of whose differentials 
respect the filtration given by the dictionary ordering on $\{0,1\}^\ell$,
and whose $d_1$ differential is obtained by adding all $\Ga_{I<I'}$
(where $I'$ is an immediate successor of $I$). 
\end{theorem}

\begin{remark}
\label{rmk:Flatten}
It is sometimes convenient to pass from ``cubically filtered''
complexes as described above to $\Z$-filtered ones. This procedure --
``flattening'' the cube (c.f.~\cite{Khovanov}) -- is done, for
example, by letting $E^1_m=
\bigoplus_{\{I\in\{0,1\}^\ell\big|\sum_{i=1}^\ell I_i=m\}}
E^1_{I}$. Note that the length of this filtration is one greater
than the
number of components of the link.
\end{remark}

\begin{remark}
Although we have stated Theorem~\ref{thm:LinkSpecSeq} for $\HFa$,
the same result can be established for $\HFp$,
with some notational changes.
\end{remark}

Before proceeding to the proof, we indulge in a purely
homological-algebraic digression. The algebraic here was inspired by a
conversation with Paul Seidel, who communicated to us some version of
Lemma~\ref{lemma:MCone}.

\subsection{Mapping cones}
\label{sec:MCone}

We begin with some terminology.

Let $A_1$ and $A_2$ be a pair of chain complexes.  A chain map
$$\phi\colon A_1 \longrightarrow A_2$$ is called a {\em
quasi-isomorphism} if the induced map on homology is an
isomorphism. Two chain complexes $A_1$ and $A_2$
are said to be {\em quasi-isomorphic}
if there is a third chain complex $B$ and 
a pair of quasi-isomorphisms $\phi_1\colon A_1 \longrightarrow B$ and 
$\phi_2\colon A_2 \longrightarrow B$.

Recall that if we have a chain map between chain complexes
$f_1\colon A_1 \longrightarrow A_2$, we can form its mapping cone
$\MCone(f_1)$, whose underlying module is the
direct sum $A_1\oplus A_2$, endowed with the differential 
$$\partial =
\left(
\begin{array}{ll}
\partial_1 & 0 \\
f_1 & \partial_2
\end{array}
\right),$$
where here $\partial_i$ denotes the differential for the chain complex
$A_i$.
(When working with coeffients in characteristic different from $2$, 
for the above
formula to work, we adopt the convention that for a chain map,
$\partial_2 \circ f_1 + f_1\circ \partial_1=0$.)
Recall that there is a short exact sequence of chain complexes
$$
\begin{CD}
0@>>> A_2@>{\iota}>> \MCone(f_1) @>{\pi}>> A_1 @>>> 0.
\end{CD}$$
This induces a long exact sequence, for which the connecting
homomorphism is the map on homology induced by $f_1$.

The mapping cylinder is natural in the following sense.
Suppose that we have a diagram of chain complexes
$$
\begin{CD}
A_1 @>{f_1}>> A_2 \\
@V{\psi_1}VV @VV{\psi_2}V \\
B_1 @>{g_1}>> B_2 \\
\end{CD}
$$
which commutes up to homotopy, then there is an induced map 
$$m(\psi_1,\psi_2)\colon \MCone(f_1) \longrightarrow \MCone(g_1)$$
which fits into the following diagram, where the rows are exact and
the squares are homotopy-commutative:
$$
\begin{CD}
0@>>> A_2@>>> \MCone(f_1) @>>> A_1 @>>> 0 \\
&& @V{\psi_2}VV @V{m(\psi_1,\psi_2)}VV @VV{\psi_1}V \\
0@>>> B_2@>>> \MCone(g_1) @>>> B_1 @>>> 0.
\end{CD}
$$

\begin{lemma}
\label{lemma:MCone}
Let $\{A_i\}_{i=1}^\infty$ be a collection of chain maps and let
$$\{f_i\colon A_i\longrightarrow A_{i+1}\}_{i\in\Z}$$
be a collection of chain maps  satisfying the following 
two properties:
\begin{enumerate}
\item
\label{item:NearlyChain}
$f_{i+1}\circ f_i$ is chain homotopically trivial, by a 
chain homotopy
$$H_i\colon A_i\longrightarrow A_{i+2}$$
\item
\label{item:Composite}
the difference
$$f_{i+2}\circ H_i-H_{i+1}\circ f_i\colon A_i \longrightarrow A_{i+3}$$ 
is a quasi-isomorphism.
\end{enumerate}
Then, $\MCone(f_2)$ is quasi-isomorphic to $A_4$.
\end{lemma}

\begin{proof}
Hypothesis~\eqref{item:NearlyChain} shows that the maps induced on
homology by the $f_i$ fit together to form a chain complex:
\begin{equation}
\label{eq:CxHom}
\begin{CD}
...H_*(A_1) @>{f_1}>> H_*(A_2)@>{f_2}>> H_*(A_3) @>{f_3}>> H_*(A_4) ...
\end{CD}
\end{equation}

Letting $\psi_i=f_{i+2}\circ H_i+H_{i+1}\circ f_i\colon A_i
\longrightarrow A_{i+3}$, then Hypothesis~\eqref{item:NearlyChain}
proves that $\psi_i$ is a chain map; and indeed that the square 
\begin{equation}
\label{eq:HoSquare}
\begin{CD}
A_i@>{f_i}>> A_{i+1} \\
@V{\psi_i}VV @VV{\psi_{i+1}}V \\
A_{i+3} @>{f_{i+3}}>> A_{i+4}
\end{CD}
\end{equation}
commutes up to homotopy. 

Next, define 
$\alpha_i\colon M(f_{i}) \longrightarrow A_{i+2}$ by:
$$\alpha_i(a_i,a_{i+1})=H_i(a_i)-f_{i+1}(a_{i+1}),$$
and 
$\beta_i\colon A_i\longrightarrow M(f_{i+1})$
by
$$\beta_i(a_i)=(f_i(a_i),H_i(a_i)).$$
Now, $$\alpha_{i+1}\circ \beta_i=\psi_i,$$
which is a quasi-isomorphism. 

Moreover, consider the diagram:
\begin{equation}
\label{eq:HoLadder}
\begin{CD}
A_2 @>{f_2}>> A_3 @>{\iota_3}>> \MCone(f_2) @>{\pi_3}>> A_2 @>{f_2}>> A_3\\
@V{=}VV @V{=}VV @V{\alpha_2}VV @V{\psi_2}VV @V{\psi_3}VV \\
A_2 @>{f_2}>> A_3 @>{f_3}>> 
A_4 @>{f_4}>> A_5 @>{f_5}>>  A_6\\
@V{\psi_2}VV @V{\psi_3}VV @V{\beta_4}VV @V{=}VV @V{=}VV \\
A_5 @>{f_5}>> A_6 @>{\iota_6}>> \MCone(f_5) @>{\pi_5}>> A_5 @>{f_5}>> 
A_6
\end{CD}
\end{equation}
The map $X\colon M(f_2)\longrightarrow A_5$ defined by
$$X(a_2,a_3)=H_2(a_2) + f_3(a_3)$$ gives a chain homotopy between
$\psi_2\circ \pi$ and $f_4\circ \alpha$, while the map $Y \colon A_3
\longrightarrow M(f_5)$ defined by $$Y(a_3)=(H_3 a_3, 0)$$ gives a
chain homotopy between $f_3\circ \beta_4$ and $\iota\circ \psi_3$.
Thus, all the squares in Diagram~\eqref{eq:HoLadder} commute up to
homotopy, and the maps induced on homology on the top and bottom rows
are exact.  From the five-lemma, it follows that the map induced on
homology $\beta_4\circ \alpha_2$ is also an isomorphism. Thus (in view
of the fact that $\alpha_{5}\circ\beta_4$ is a quasi-isomorphism, we
conclude $\beta_4$ and hence $\alpha_2$ is a quasi-isomorphism.
\end{proof}

It is useful to interpret Lemma~\ref{lemma:MCone} in the following
terms. Under the hypotheses of that lemma,
we can form an ``iterated mapping cone'' 
$\MCone(f_1,f_2,f_3)$ whose underlying module is
$A_1\oplus A_2 \oplus A_3$, and whose differential is given by the matrix
\begin{equation}
\label{eq:Reformulate}
\partial=
\left(
\begin{array}{lll}
\partial_1 & 0 & 0 \\
f_1 & \partial_2 & 0 \\
H_1 & f_2 & \partial_3
\end{array}
\right).
\end{equation}
Indeed, Hypothesis~\eqref{item:NearlyChain} guarantees that 
$\partial$ determines a differential on
$\MCone(f_1,f_2,f_3)$.  
Lemma~\ref{lemma:MCone} can be interpreted as saying that
$H_*(\MCone(f_1,f_2,f_3))=0$. Now, there is a short exact sequence
$$
\begin{CD}
0@>>> A_3@>>> \MCone(f_1,f_2,f_3) @>>>\MCone(f_1) @>>>0,
\end{CD}$$
and the connecting homomorphism $$H_*(\MCone(f_1))\longrightarrow
H_*(A_3)$$ is the map on homology induced by the map $\alpha_1$ in the
above lemma.

\subsection{Pseudo-holomorphic $n$-gons}
\label{sec:NGons}

With the above homological algebra in place, we proceed to the geometrical
underpinnings of Theorem~\ref{thm:LinkSpecSeq}.
We will make heavy use of the pseudo-holomorphic polygon construction,
c.f.~\cite{FOOO}, \cite{QAssoc}, \cite{Silva} and its relationship
with Heegaard Floer homology, as explained in
Section~\ref{HolDisk:sec:HolTriangles} of~\cite{HolDisk}. We recall
this construction very briefly here.

Let $\Sigma$ be a connected, closed, oriented two-manifold of genus $g$,
and fix $(m+1)$ $g$-tuples of attaching circles $\{\etas^{i}\}_{i=1}^m$.
Specifically, for each fixed $i$, the set
$\etas^i=\{\eta^i_j\}_{j=1}^g$ is a collection of $g$
pairwise disjoint,
homologically linearly independent, embedded curves in $\Sigma$. We choose
also a reference point $\Sigma$ disjoint from all the $\eta^i_j$. 
In the terminology of~\cite{HolDisk}, 
this data $(\Sigma,\etas^0,...\etas^m,z)$ is called a pointed Heegaard
$(m+1)$-tuple.
We have a map of groups 
$$\fa_{\eta^0,...,\eta^m}\colon
\bigotimes_{i=1}^m\CFa(Y_{\eta^{i-1},\eta^i}) \longrightarrow
\CFa(Y_{\eta^0,\eta^m}),$$
where here $Y_{\eta^i,\eta^j}$ denotes the
three-manifold described by the Heegaard diagram
$(\Sigma,\etas^i,\etas^j)$.  This map is obtained by counting
pseudo-holomorphic $m+1$-gons in $\Sym^g(\Sigma)$ which are disjoint
from the subvariety $\{z\}\times\Sym^{g-1}(\Sigma)$. 

More precisely, let ${\mathbb T}_{\eta^i}$ denote the $g$-dimensional
torus $\eta^i_1\times...\times\eta^i_g$ in the $g$-fold symmetric
product $\Sym^g(\Sigma)$. A {\em Whitney $(m+1)$-gon} is a map
$u$ from the standard $(m+1)$-gon into $\Sym^g(\Sigma)$
which maps the $i^{th}$ edge into ${\mathbb T}_{\eta^i}$ (where here the edges are labelled $0,...,m$). Fixing 
$\x_i\in{\mathbb T}_{\eta^{i-1}}\cap{\mathbb T}_{\eta^{i}}$ and
$\y\in{\mathbb T}_{\eta^0}\cap{\mathbb T}_{\eta^{m}}$,
we let $\pi_2(\x_1,...,\x_m,\y)$ denote the set of homotopy classes
of Whitney $(m+1)$-gons
which, 
for $i=1,...m$, map the vertex between the $(i-1)^{st}$ and the $i^{th}$
edge to $\x_i$ and the vertex between
$m^{th}$ and $0^{th}$ edges to $\y$.

For fixed $\varphi\in\pi_2(\x_1,...,\x_m,\y)$, we let $\ModFlow(\varphi)$
denote the set of pseudo-holomorphic representatives for $\varphi$.
With this notation in place, then, the map $\fa_{\eta^0,...,\eta^m}$ (when $m>1$) is defined by
$$\fa_{\eta^0,...,\eta^m}(\x_1\otimes...\otimes\x_m)=\sum_{\y\in{\mathbb T}_{\eta^0}\cap{\mathbb
    T}_{\eta^m}}
\sum_{\{\varphi\in\pi_2(\x_1,...,\x_m,\y)\big|\Mas(\varphi)=0,n_z(\varphi)=0\}}
\big(\#\ModFlow(\varphi)\big)\cm\y,$$
where here $n_z(\varphi)$ denotes the intersection number of $\varphi$
with the subvariety $\{z\}\times \Sym^{g-1}(\Sigma)\subset
\Sym^g(\Sigma)$, and $\Mas(\varphi)$ denotes the expected dimension
of the moduli space $\ModFlow(\varphi)$ (i.e.  the Maslov index of
$\varphi$, c.f.~\cite{FloerMaslov}, \cite{SalamonZehnder}).  In the
special case where $m=1$, we sum over homotopy classes wih
$\Mas(\varphi)=1$, and count points in the quotient space
$\ModFlow(\varphi)/\R$. Thus, when $m=1$, the map
$\fa_{\eta^0,\eta^1}$ is simply the differential for the Heegaard Floer
chain complex for $Y_{\eta^0,\eta^1}$, and when $m=2$,
$\fa_{\eta^0,\eta^1,\eta^2}$ is the chain map induced from the counts of
pseudo-holomorphic triangles.

These maps are well-known to satisfy a generalized associativity
property, c.f.~\cite{FOOO}, \cite{QAssoc}, \cite{Silva}:
\begin{equation}
\label{eq:Associativity}
\sum_{0\leq i < j \leq m}
\fa_{\eta^{i},\eta^{i+1},...\eta^{j}}\circ
\fa_{\eta^0,...\eta^{i-1},\eta^i,\eta^j,...,\eta^m} = 0.
\end{equation}
For example, when $m=1$, the above associativity statement is the is
equivalent to the fact that the square of the differential for
$\CFa(Y_{\eta^0 \eta^1})$
is trivial.
When $m=2$, associativity asserts that the maps induced by holomorphic
triangles are chain maps, and when $m=3$, it states that the triangle
pairing is associative, up to chain homotopy (hence the name).
Thinking of the tuples $\etas^0,...,\etas^m$
as corresponding to edges of an $(m+1)$-gon, we see that
for any pair of edges, there is
a corresponding degeneration of the $(m+1)$-gon as a juxtaposition of a 
pair of $a$- and $b$-gons, with $a+b-3=m$. The above sum a composition
over all such degenerations of the corresponding composition of maps.

We can construct Heegaard diagrams for the $Y(I)$ as follows.  Given
$Y$ with the framed link $L$, we can construct a Heegaard triple
$(\Sigma,\alphas,\betas,\gammas,z)$, where here
$\beta_{1},...,\beta_\ell$ are meridians for the links,
$\gamma_1,...,\gamma_\ell$ are corresponding framing curves (and
$\gamma_{\ell+1},...,\gamma_g$ are exact Hamiltonian translates of
$\beta_{\ell+1},...,\beta_g$).  (For more on the construction this
diagram, see Section~\ref{HolDiskFour:sec:ECobord}
of~\cite{HolDiskFour}. In the terminology of that paper, the Heegaard
triple we are considering is the Heegaard triple subordinate to some
bouquet for the framed link $L$.)  We also choose curves
$\delta_1,...,\delta_\ell$ to be corresponding curves representing the
framings obtained by adding meridians to the original framings (and
$\delta_{\ell+1},...,\delta_{g}$ are exact Hamiltonian translates of
the $\beta_{\ell+1},...,\delta_{g}$).  Given
$I\in\{0,1,\infty\}^\ell$, let $\etas(I)=\{\eta_1,...,\eta_g\}$ denote
the $g$-tuple of attaching circles, where here
$$\eta_i=\left\{\begin{array}{ll} \beta_i & {\text{if $i>n$ or
$m_i=\infty$}} \\ \gamma_i & {\text{if $m_i=0$}} \\ \delta_i &
{\text{if $m_i=1$.}}
\end{array}\right.$$
Thus, a Heegaard diagram for $Y(\MulFram)$ is given by
$(\Sigma,\alphas,\etas(\MulFram),z)$.

Given a sequence of multi-framings $\MulFram^0<...<\MulFram^k$,
there is an induced map
$$D_{\MulFram^0<...<\MulFram^k}
\colon \CFa(Y(I^0))\longrightarrow \CFa(Y(I^k))$$
defined by 
\begin{equation}
\label{eq:DefDiff}
D_{\MulFram^0<...<\MulFram^k}(\xi) =
\fa_{\alpha,\eta(\MulFram^0),...,\eta(\MulFram^k)}(\xi\otimes
\Thetaa_1\otimes ...\otimes \Thetaa_k),
\end{equation}
where $\Thetaa_i$ are cycles
representing the canonical top-dimensional
generators for $\HFa$ of
$Y_{\eta(I^i),\eta(I^{i+1})}$, which is a connected
sum of several copies of $(S^2\times S^1)$.
(We will be lax about distinguishing here between 
intersection points
in ${\mathbb T}_{\eta(I^i)}\cap {\mathbb T}_{\eta(I^{i+1})}$
and generators of the homology groups
$\HFa(\#^k(S^2\times S^1))$; after all,
one can equip this three-manifold with a Heegaard diagram for which
the intersection points all correspond to closed chains, and hence each
generator is represented by a unique intersection point.)

Let $X=\bigoplus_{I\in \{0,1,\infty\}^\ell} \CFa(Y(I))$, endowed with the map
$$D \colon X \longrightarrow X,$$ defined by $$D \xi =
\sum_{J}\sum_{\{I=I^1<...<I^k=J\}} D_{I^1<...<I^k}(\xi),$$
where here  the
index set of the inner sum is the set of all increasing sequences
connecting $I$ to $J$, with the property that for all $i=1,...,k+1$,
$I^{i+1}$ is an immediate successor of $I^{i}$.

\begin{lemma}
\label{lemma:ModelCalc}
Fix $I, J\in \{0,1,\infty\}^\ell$.
We have that
$$\sum_{I=I^0<I^1<...<I^k=J}\fa_{\eta(I^0),...,\eta(I^k)}(\Thetaa_1\otimes...\otimes\Thetaa_k)\equiv 0,$$
where again the sum is taken over sequences with the property that $I^{i+1}$
is an immediate successor of $I^i$.
\end{lemma}

\begin{proof}
We consider the case where $k>2$. In this case, there is a juxtaposition
of triangles representing
$$\fa_{\eta(I^0),\eta(I^1)\eta(I^2)}\circ 
\fa_{\eta(I^0),\eta(I^2),\eta(I^3)}\circ ...
\circ \fa_{\eta(I^0),\eta(I^i),\eta(I^{i+1})}\circ ... 
\circ \fa_{\eta(I^0),\eta(I^{k-1}),\eta(I^{k})}.
$$ This juxtaposition gives rise to an $(k+1)$-gon 
$\varphi\in\pi_2(\Thetaa_1,...,\Thetaa_k,\Theta)$,
with $\cald(\varphi)\geq 0$ and $n_z(\varphi)=0$, 
where here $\Theta\in{\mathbb T}_{\eta(I^0)}\cap {\mathbb
T}_{\eta(I^k)}$.  By additivity of the Maslov index, this $k+1$-gon has
$\mu(\varphi)=k-2$; it is also easily seen to have
$n_z(\phi)=0$. It follows  that there are no $k+1$-gons
$\varphi'\in \pi_2(\Thetaa_1,...,\Thetaa_k,\Thetaa')$ with
$\Mas(\varphi')=0$ and $\cald(\phi')\geq 0$. 

Consider now the case where $k = 2$. In this case, $I$ and $J$ differ
in at least one place, and at most two.

If $I$ and $J$ differ in one place, a direct inspection of the
Heegaard triple (which leads to the ``blowup formula'' in
~\cite{HolDiskFour}) shows that the maps cancel in pairs.  This is
spelled out in Proposition~\ref{HolDiskTwo:prop:HoClassesCancel}
of~\cite{HolDiskTwo}.

If they differ in two places, there are two choices for $I^1$ with
$I=I^0<I^1<J=I^2$. For each possible $I^1$, it is the case that
$$\fa_{\eta(I^0),\eta(I^1),\eta(I^2)}(\Thetaa_1\otimes\Thetaa_2)=\pm \Thetaa_3.$$
One can see this by explicitly drawing the Heegaard triple, which
splits into torus summands, as in~\cite{HolDiskFour}. (See also
Proposition~\ref{prop:STwoSOne} below.)  Since there are two choices,
they cancel in the sum, modulo two. 

Indeed, it is straightforward to give all triangles orientations so
that the cancellation takes place with coefficients in $\Z$, as
well. For example, for fixed $I$ and $J$ which differ in two places,
and for which we have a sequence of immediate successors of the form
$I<K<J$; we have exactly one other sequence of immediate successors
$I<L<J$. Label these two possible intermediate mult-framings so that
$K<L$. We then define orientations on the triangles (and choose our canonical
generators for all the $Y_{\eta(I),\eta(I')}$) so that
$$\fa_{\eta(I),\eta(K),\eta(J)}(\Thetaa_1\otimes\Thetaa_2) =\Thetaa_3,$$
while $$\fa_{\eta(I),\eta(L),\eta(J)}(\Thetaa_1\otimes\Thetaa_2)
=-\Thetaa_3.$$

In the case where $k=1$, the stated relation is simply the one that
$\Thetaa_1$ is a cycle.
\end{proof}

\begin{prop}
\label{prop:DisDiff}
The  map $D$ from Equation~\eqref{eq:DefDiff}
satisfies $D^2=0$.
\end{prop}

\begin{proof}
This follows from the associativity formula (Equation~\eqref{eq:Associativity})
for the
Heegaard tuple
$(\Sigma,\alphas,\etas(I^1), ...., \etas(I^k),z)$,
together with Lemma~\ref{lemma:ModelCalc}. Specifically, 
according to that lemma, 
the only degenerations in Equation~\eqref{eq:Associativity}
which do not contribute $0$ to the sum
are the ones which involve $\alphas$ in both polygons.
Those, in turn, are the various components of $D^2$.

Strictly speaking, for this result to hold over $\Z$ (not just over
$\Zmod{2}$), note that Lemma~\ref{lemma:ModelCalc} holds for
particular orientations on the spaces of triangles for
$$(\Sigma,\eta(I^0),\eta(I^1),\eta(I^2)),$$
where here $I^0<I^1<I^2$ are immediate
successors. In fact, it is straightforward to see that such a system
of orientations can be used to induce a coherent orientation on all
the polygons, so that associativity still holds, and hence the result holds
over $\Z$, as well.
\end{proof}

In view of Proposition~\ref{prop:DisDiff}, we can think of
$X$ is a chain complex, endowed with the differential $D$.
We can define some other associated complexes as follows.
If $S\subset \{0,1,\infty\}^\ell$ is a subset
with the property that for each $I,J\in S$, for
all $K\in\{0,1,\infty\}^\ell$ with $I<K<J$, we also have that $K\in S$,
then we let $X(S)$ denote the group 
$\bigoplus_{I\in S}\CFa(Y(I))$ endowed with the differential naturally induced by $D$.

With all the notational background, we are now ready to prove a strong form
of the surgery long exact sequence for a single knot in a three-manifold $Y$.

\begin{theorem}
\label{thm:KnotCase}
Let $K$ be a framed knot in a three-manifold $Y$, and let 
$$\fa\colon \CFa(Y_0(K))\longrightarrow \CFa(Y_1(K))$$
denote the chain map
induced by the cobordism. Then, the chain complex $\CFa(Y)$ is
quasi-isomorphic to the
mapping cone of $\fa$.
\end{theorem}

\begin{proof}
  To start, let $(\Sigma,\alphas,\betas,\gammas,\deltas,z)$ denote the
  associated Heegaard quintuple. In particular, $Y_{\alpha,\beta}$,
  $Y_{\alpha,\gamma}$ $Y_{\alpha,\delta}$ describe $Y$, $Y_0$, and
  $Y_1$ respectively, and the remaining three-manifolds on the
  boundary describe $\#^{g-1}(S^2\times S^1)$.  Indeed, to fit
  precisely with the hypotheses of that lemma, we choose infinitely
  many copies of the $g$-tuples $\betas$, $\gammas$, and $\deltas$
  (denoted $\betas^{(i)}$, $\gammas^{(i)}$, $\deltas^{(i)}$ for
  $i\in\Z$), all of which are generic exact Hamiltonian perturbations
  of one another, in the interest of admissibility (in the sense of
  Subsection~\ref{HolDisk:subsec:Admissibility} of~\cite{HolDisk}).
  
  In this case, the chain map we described earlier $X$ splits (as a
  module) as $\CFa(Y_{0})\oplus \CFa(Y_1)
  \oplus\CFa(Y)$, and its
  differential decomposes as \begin{equation} \label{eq:BigDiff}
  \partial= \left(\begin{array}{lll} D_{0} & 0 & 0 \\ D_{0<1} & D_{1}
  & 0 \\ D_{0<1<\infty} & D_{1<\infty} & D_{\infty} \end{array} \right).
  \end{equation}
  
  Letting $\CFa(Y)$, $\CFa(Y_0)$, and $\CFa(Y_1)$ play the roles of
  $A_1$, $A_2$, and $A_3$ respectively, the various components
  of the differential play the roles of the $f_i$ and $H_i$
  (compare Equations~\eqref{eq:BigDiff} and~\eqref{eq:Reformulate}),

  Indeed, $A_{3i+1}$, $A_{3i+2}$ and $A_{3i+3}$ all represent
  $\CFa(Y_0)$, $\CFa(Y_1)$ and $\CFa(Y)$ respectively, only now we use the
  various translates of the $\betas$, $\gammas$, and $\deltas$;
  in particular $A_{3i+1}$ is the Floer complex $\CFa(\alphas,\betas^{(i)})$.
  
  Hypothesis~\eqref{item:NearlyChain} of Lemma~\ref{lemma:MCone}
  follows at once from the fact that $D$ is a chain complex
  (Proposition~\ref{prop:DisDiff}).

  It remains to verify Hypothesis~\eqref{item:Composite} of
  Lemma~\ref{lemma:MCone}.

  Let $\theta_i$ be the chain homotopy equivalences induced by
  equivalences of Heegaard diagrams; e.g.  $\theta_{3i}$
  is the chain map $\CFa(\alphas,\betas^{(i)})
  \longrightarrow\CFa(\alphas,\betas^{(i+1)})$ 
  obtained by product with the canonical generator
  $\Thetaa_{\beta^{(i)},\beta^{(i+1)}}$). 

  We claim that $$f_{3}\circ H_1+H_{2}\circ f_1\colon
  A_1\longrightarrow A_4$$ is chain homotopic to $\theta_1$, and the
  chain homotopy is given by $$ \x\mapsto
  f_{\alpha,\beta,\gamma,\delta,\beta^{(i)}}(\x\otimes
  \Thetaa_{\beta,\gamma}\otimes\Thetaa_{\gamma\delta}\otimes
  \Thetaa_{\delta,\beta^{(1)}}).  $$ This in turn follows at once from
  associativity, together with the fact that 
  \begin{equation}
        \label{eq:CanonicalSquare}
  f_{\beta,\gamma,\delta,\beta^{(1)}} (\Thetaa_{\beta,\gamma}\otimes
  \Thetaa_{\gamma,\delta}\otimes \Thetaa_{\delta,\beta^{(1)}})=
  \Thetaa_{\beta,\beta^{(1)}}.  
  \end{equation}
   This latter equality follows from a
  direct inspection of the Heegaard diagram for the quadruple
  $(\Sigma,\betas,\gammas,\deltas,\betas^{(1)},z)$.  (i.e. the count
  of pseudo-holomorphic pentagons), as illustrated in
  Figures~\ref{fig:Square} and ~\ref{fig:RestOfSquare}.

\begin{figure}
\mbox{\vbox{\epsfbox{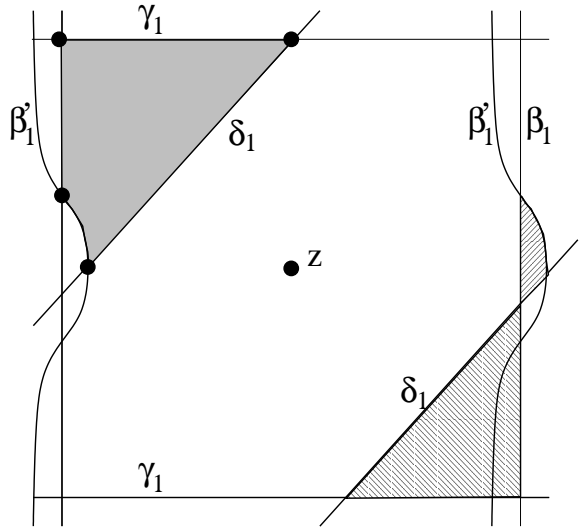}}}
\caption{\label{fig:Square}
{\bf A holomorphic quadrilateral.}
The shaded quadrilateral has a unique holomorphic representative
(by the Riemann mapping theorem), while the one indicated with
the  hatching does not,
as it has both positive and negative local multiplicities,
as indicated by the two directions in the hatching.}
\end{figure}


\begin{figure}
\mbox{\vbox{\epsfbox{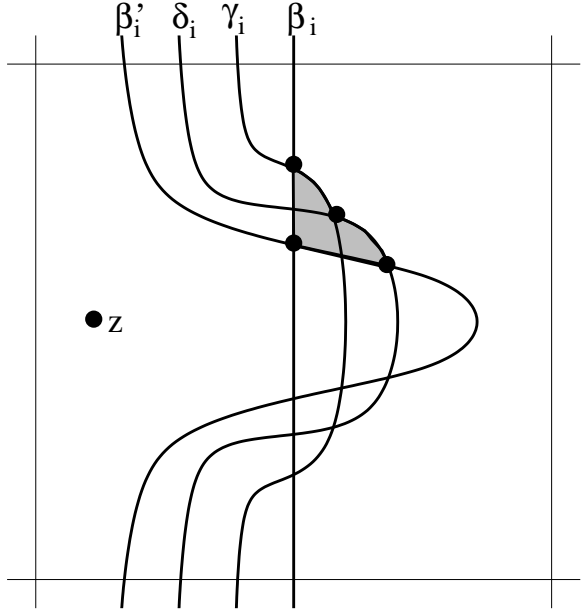}}}
\caption{\label{fig:RestOfSquare}
  {\bf Other factors of the holomorphic quadrilateral.}  We have
  illustrated here a Heegaard quadruple (in a genus one surface) whose
  four boundary components are $S^2\times S^1$. In the homotopy class
  indicated by the shaded quadrilateral $\varphi_i\in
  \pi_2(\Thetaa_{\beta_i,\gamma_i},\Thetaa_{\gamma_i,\delta_i},\Thetaa_{\delta_i,\beta_i'},\Thetaa_{\beta_i',\beta_i})$,
  there is a moduli space of pseudo-holomorphic quadrilaterals which
  is clearly one-dimensional, parameterized by a cut at the vertex
  where $\gamma_i$ and $\delta_i$ meet.  We take the connected sum of
  $g-1$ copies of this picture (at the reference point $z$) with the
  picture illustrated in Figure~\ref{fig:Square} to obtain the general
  case of the quadrilateral considered in the proof of
  Theorem~\ref{thm:KnotCase}.}
\end{figure}

  In Figure~\ref{fig:Square}, we consider the special case where the
  genus $g=1$. In the picture, and in the following discussion,
  $\beta^{(1)}_1$ is denoted $\beta'_1$. The four corners of the
  shaded quadrilateral are the canonical generators
  $\Thetaa_{\beta_1,\gamma_1}$, $\Thetaa_{\gamma_1,\delta_1}$,
  $\Thetaa_{\delta_1,\beta'_1}$, and $\Thetaa_{\beta'_1,\beta_1}$ (read
  in clockwise order). Indeed, it is straightforward to see (by
  passing to the universal cover), that the shaded quadrilateral
  represents the only homotopy class $\varphi_1$ of Whitney
  quadrilaterals with $n_z(\varphi_1)=0$ and all of whose local
  multiplicities are non-negative.  By the Riemann mapping theorem,
  now, this homotopy class $\varphi_1$ has a unique holomorphic
  representative $u_1$. (By contrast, we have also pictured here 
  another Whitney quadrilateral with  hatchings, 
  whose local multiplicities are all
  $0$, $+1$, and $-1$; $+1$ at the region where the hatchings go in
  one direction and $-1$ where they go in the other.)

  For the general case ($g>1$), we take the connected sum of the case
  illustrated in Figure~\ref{fig:Square} with $g-1$ copies of the
  torus illustrated in Figure~\ref{fig:RestOfSquare}. In this picture,
  we have illustrated the four curves $\beta_i$, $\gamma_i$,
  $\delta_i$, $\beta'_i$ for $i>1$, which are Hamiltonian translates
  of one another.  Now, there is a homotopy class of quadrilateral
  $\varphi_i\in
  \pi_2(\Thetaa_{\beta_i,\gamma_i},\Thetaa_{\gamma_i,\delta_i},
  \Thetaa_{\delta_i,\beta_i'},\Thetaa_{\beta_i',\beta_i})$, and a
  forgetful map $\ModFlow(\phi)\longrightarrow \ModRect$ which
  remembers only the conformal class of the domain (where here
  $\ModRect$ denotes the module space of conformal classes of disks
  with four marked boundary points, also referred to simply as
  quadrilaterals). Both moduli spaces are one-dimensional (the first
  moduli space is parameterized by
  the length of the cut into the region, while the second is
  parameterized by the ratio of the length to the width, after the
  quadrilateral is uniformized to a rectangle). By Gromov's
  compactness theorem, the forgetful map is proper; and it is easy to
  see that this it has degree one, and hence for some generic
  conformal class of quadrilateral, there is a unique
  pseudo-holomorphic quadrilateral whose domain has the specified
  conformal class. Now, letting $u_1$ (the pseudo-holomorphic
  representative of the homotopy class $\varphi_1$ described in the
  previous paragraph) determine the conformal class of the rectangle,
  we let $u_i$ for $i>1$ be the pseudo-holomorphic representatives for
  $\phi_i$ whose domain supports the same conformal class.  Then
  $u_1\times... \times u_g\in \varphi_1\times...\times\varphi_g$ is
  easily seen to be the unique holomorphic quadrilateral in
  $\pi_2(\Thetaa_{\beta,\gamma},\Thetaa_{\gamma,\delta},
  \Thetaa_{\delta,\beta'},\Thetaa_{\beta',\beta})$, hence proving
  Equation~\eqref{eq:CanonicalSquare} which, in turn, yields
  Hypothesis~\eqref{item:Composite} of Lemma~\ref{lemma:MCone}.  The
  theorem now follows directly from Lemma~\ref{lemma:MCone}.
\end{proof}

We now turn to Theorem~\ref{thm:LinkSpecSeq}.

\vskip.2cm
\noindent{\bf Proof of Theorem~\ref{thm:LinkSpecSeq}.}
The theorem  is established by induction on the number of components
of the link. The case where the link has a single component is
a direct consequence of Theorem~\ref{thm:KnotCase}.

We form the chain complex $X$ as before. We claim first that
$H_*(X)=0$. 

Let $S$ be any subset of $\{0,1,\infty\}^{\ell}$ with the following
property. For each $\sigma\in S$, there are two other elements $\tau,
\nu\in S$, with the property that the three subsets $\sigma,\tau,\nu$
are three distinct elements, and they all differ in a single
component. It is a straightforward application of
Theorem~\ref{thm:KnotCase} (or, more precisely, the reformulation of
the mapping cone lemma, Lemma~\ref{lemma:MCone}, described after its
proof) that $H_*(X(S))=0$.  In particular, $H_*(X)=0$.

In particular, we can let 
$$S=\{0,1\}^{\ell-1}\times \{0,1,\infty\}.$$
By the above criterion, $H_*(X(S))=0$. Indeed, the 
letting $T=\{0,1\}^\ell$, we have a short exact sequence
$$
\begin{CD}
0@>>> X(T) @>{f}>> X(S) @>{g}>> X(\{0,1\}^{\ell-1}\times\{\infty\})@>>>0
\end{CD}$$
from which it follows at once that the connecting homomorphism induces
an isomorphism in homology
$$H_*(X(\{0,1\}^{\ell-1}\times\{\infty\}))\stackrel{\cong}\longrightarrow
H_*(X(\{0,1\}^\ell)).$$ By our inductive hypothesis, it follows that
$H_*(X\{0,1\}^{\ell-1}\times\{\infty\})\cong\HFa(Y)$, completing the proof.
\qed

\section{Khovanov's invariants}
\label{sec:Khovanov}

We briefly describe here Khovanov's categorification of the Jones
polynomial; for more details, see~\cite{Khovanov} and~\cite{BarNatan}.
We make some simplifying assumptions here: we will use coefficients in
$\Zmod{2}$ throughout, and we specialize to the case where $U=0$ (in
Khovanov's sense). Our notation and exposition are taylored to fit
neatly into the context of the present paper. In particular, the groups
we describe here are actually the Khovanov homology of the mirror of $L$.

Let $X=S_1\cup ...\cup S_k$ be a collection of disjoint
embedded, simple closed curves in the plane. Let $Z(X)$ denote the
$\Zmod{2}$-vector space, formally generated by the components
$[S_1],...,[S_k]$, and let $V(X)$ denote
exterior algebra $$V(X)=\Wedge^* Z(X),$$ 
i.e. this is the quotient of the polynomial algebra over $\Zmod{2}$
generated by $[S_i]$, divided out by the relations $[S_i]^2=0$ for $i=1,...,k$.
(When comparing the notation used in the 
discussion here with that of~\cite{Khovanov},
observe that the element $[S_{1}]\wedge...\wedge [S_{m_\ell}]\in V(X)$, 
where $\{m_j\}_{j=1}^\ell$
is a subsequence of $\{1,...,k\}$ corresponds to the element
$$v_{1}^{\epsilon_1}\otimes...\otimes v_k^{\epsilon_k}$$
where here $\epsilon_i\in \{\pm\}$ is obtained by
$$
\epsilon_i=
\left\{\begin{array}{ll}
- & {\text{if $i\in\{m_j\}_{j=1}^\ell$}} \\
+ & {\text{otherwise}}
\end{array}
\right.
$$
in Khovanov's notation, c.f.~\cite{Khovanov}.)

Next, consider a pair of pants, thought of as a morphism from
$X=S_1\cup...\cup S_k\cup S_{k+1}$ to a new submanifold
$X'=S_1\cup...\cup S_{k-1}
\cup S_k'$ containing a component $S_k'$  which 
is obtained by merging $S_k$ and $S_{k+1}$. In this case, we have a
natural identification $$Z(X')=Z(X)/[S_{k}]\sim [S_{k+1}],$$ and
correspondingly natural isomorphisms
\begin{eqnarray*}
\alpha\colon (S_{k+1}-S_k)\wedge V(X)\stackrel{\cong}{\longrightarrow}
V(X')
&{\text{and}}&
\beta\colon V(X')\stackrel{\cong}{\longrightarrow}
V(X)/(S_{k+1}-S_k)\wedge V(X).
\end{eqnarray*}
We then define the multiplication $$m\colon
V(X)\longrightarrow V(X')$$ to be the
composite $$\begin{CD} V(X) @>{(S_{k+1}-S_k)\wedge\cdot}>>
(S_{k+1}-S_{k})\wedge V(X) @>{\alpha}>> V(X').
\end{CD}
$$

By reversing the orientation of the ``pair of pants'', we have a
morphism from $X'$ to $X$, instead.  In this case, we have a
comultiplication $$\Delta \colon V(X') \longrightarrow V(X)$$ induced
by the composition $$\begin{CD} V(X') @>{\beta}>> 
\frac{V(X)}{(S_{k+1}-S_{k})\wedge
V(X)}@>{(S_{k+1}-S_k)\wedge}>> V(X).
\end{CD}
$$

Let $L$ be a link, and fix a generic projection of $L$, $\Diagram$,
with ${\ell}$ double points.  One can form resolutions indexed by
subsets $I\in
\{0,1\}^\ell$.  Specifically, for each $I$, $\Diagram(I)$ is a
disjoint union of circles in the plane.  If $I'$ is an immediate
successor of $I$, then $\Diagram(I')$ differs from $\Diagram(I)$ by a
single pair of pants. 
We have a map
$$\dkh_{I<I'}\colon V(\Diagram(I))
\longrightarrow V(\Diagram(I')),$$
given by multiplication or co-multiplication, according to whether
$\Diagram(I')$ has one fewer or one more component than $\Diagram(I)$.

\begin{figure}
\mbox{\vbox{\epsfbox{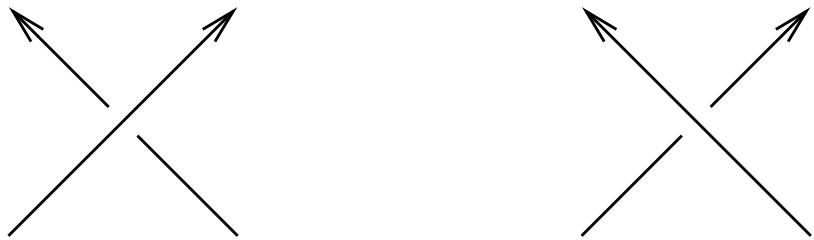}}}
\caption{\label{fig:CrossingSigns}
{\bf{Crossing conventions.}} Crossings of the first kind are 
positive, and those of the second kind are negative.}
\end{figure}

Fix a diagram $\Diagram$ for an oriented link $L$, and let
$n_+(\Diagram)$ resp. $n_-(\Diagram)$ denote the number of 
positive resp. negative crossings for the link $L$, according to the 
usual conventions (c.f. Figure~\ref{fig:CrossingSigns}). Consider
next the graded
Abelian group $$\CKh(\Diagram,m) = \bigoplus_{\{I\in\{0,1\}^\ell
\big|\sum_{i=1}^\ell I_i+n_+(\Diagram)=m\}}V(\Diagram(I)),$$ endowed with the
differential $$\dkh\colon \CKh(\Diagram,m)\longrightarrow
\CKh(\Diagram,m-1)$$ whose restriction to $V(\Diagram(I))\subset
\CKh(L,m)$ is the sum $$\dkh=\sum_{I<I'}\dkh_{I<I'},$$ where the sum
is taken over all immediate successors $I'$ of $I$.  In each dimension
$m$, $\CKh(\Diagram,m)$ is endowed with an additional grading, the
``q-grading'', defined by the splitting
$$V(\Diagram(I))=\bigoplus_{n\in\Z}V_n(\Diagram(I)),$$
where
$$V_n(\Diagram(I))=\Wedge^{k}Z^*(\Diagram(I)),$$
where
$$n=\dim Z(\Diagram)-2k-n_-(\Diagram)+2n_+(\Diagram).$$
Correspondingly, we write
$$\CKh(\Diagram)=
\bigoplus_{m,n\in\Z}\CKh(\Diagram,m,n).$$ Note that $\dkh$ carries $\CKh(\Diagram,m,n)$
to $\CKh(\Diagram,m-1,n-1)$.  

It is easy to see that $\dkh^2=0$. 
Khovanov's homology of the mirror of $L$ is the Abelian group 
$$\Kh(r(L))=H_*(\CKh_*(\Diagram),\dkh),$$
thought of as bi-graded Abelian group
$$
\Kh(r(L))=\bigoplus_{m,n\in\Z\oplus\Z} \Kh(\Diagram,m,n).$$

Note that the complex $\CKh(\Diagram)$ depends on the projection of
$L$. Khovanov shows, however, that the homology of this complex is 
independent of this choice, i.e. $\Kh(L)$ is a link invariant.
Moreover, he shows that these groups satisfy a
skein exact sequence
$$
\begin{CD}
..@>>>\Kh(r(L))@>>>\Kh(r(L_0)) @>>>\Kh(r(L_1))@>>> ...
\end{CD}
$$

Khovanov's theory is related to the Jones polynomial by the formula
$${\widehat J}(r(L))=\sum_{m,n} (-1)^m  \Big(\rk~\Kh(L,m,n)\Big)\cm q^{n},$$
where here ${\widehat J}(L)\in\Z[q,q^{-1}]$ is the unnormalized
Jones polynomial of the link $K$, characterized by the formulas 
\begin{eqnarray*}
{\widehat J}(\emptyset)&=& 1 \\
{\widehat J}\Big({\text{(unknot)}}\cup L\Big) &
\dot{=}& (q+q^{-1})\cm {\widehat J}(L) \\
{\widehat J}(r(L)) &\dot{=}& {\widehat J}(r(L_0)) - q \cm {\widehat J}(r(L_1)),
\end{eqnarray*}
where in the last equation, $L_0$ and $L_1$ are is taken with respect 
to the two resolutions at any double-point of any projection of $L$,
and where for $f,g\in \Z[q,q^{-1}]$, we write $f\dot{=}g$ if
$f=q^j\cdot g$ for some $j\in\Z$. 

In~\cite{KhovTwo}, Khovanov gives a  modification of the above
constructions to define a ``reduced'' theory $\Khred(L)$, which is
related to the normalized Jones polynomial $J(L)$ is defined by
$(q+q^{-1})\cm J(L)={\widehat J}(L)$. 
For the reduced theory, one marks a generic point in the projection of
$L$, so that now in all the various resolutions, there is always a
distinguished circle. The reduced Khovanov complex is is the
quotient of $\CKhred(\Diagram)$
by the subcomplex of $\CKh(\Diagram)$ given by
$$\bigoplus_{I\in\{0,1\}^\ell} [S_I]\wedge V(\cald(I))$$
where here $S_I$ is the component in $\Diagram(I)$ which contains the
marked point. 
This gives a chain complex whose
which splits into summands indexed by $I\in\{0,1\}^\ell$, and the corresponding
summand is denoted
$${\widetilde V}(\Diagram(I))=V(\Diagram(I))/[S_I]\wedge V(\Diagram(I)).$$
For this construction.  $$J(r(L))=\sum_{m,n}
(-1)^m\Big(\rk~\Khred(L,m,n)\Big)\cm q^{n},$$
c.f.~\cite{KhovTwo}.

We have described Khovanov's construction with $\Zmod{2}$
coefficents. In fact, Khovanov's original definition
from~\cite{Khovanov} make sense with coefficients in $\Z$. In this
case, however the reduced homology (described in~\cite{KhovTwo})
depends on the link, together with the distinguished component
containing the marked point. If one takes coefficients in $\Zmod{2}$,
as we have here, it is easy to see that the reduced theory is
independent of this additional choice, hence giving a link invariant.

\section{The spectral sequence for a branched double cover}
\label{sec:DCovSpecSeq}

Throughout this section, we fix our coefficient ring to be $\Zmod{2}$:
i.e.  if $Y$ is a three-manifold, $\HFa(Y)$ will denote $\HFa$ of $Y$
with coefficients in $\Zmod{2}$, and similarly, $H_*(Y)$ will denote
singular homology with coefficients in $\Zmod{2}$.

In comparing Khovanov's homology with $\HFa$, we rely on the following
(fairly straightforward)
result about $\HFa$, proved in earlier papers. For the statement, note
that if $Y$ is any three-manifold, then $\HFa(Y)$ is a module over the
algebra $\Wedge^* H_1(Y)/\Tors$.

\begin{prop}
\label{prop:STwoSOne}
Let $Y\cong \#^k(S^2\times S^1)$. Then, $\HFa(Y)$ is a rank one, free
module over the ring $\Wedge^* H_1(Y)$, generated by some class
$\Theta\in\HFa(Y)$.  Moreover, if $K\subset Y$ is a curve which is
dual to one of the circles in one of the $S^2\times S^1$ summands,
then the three-manifold $Y'=Y_0(K)$ is diffeomorphic to
$\#^{k-1}(S^2\times S^1)$, with a natural identification $$\pi\colon
H_1(Y)/[K]\longrightarrow H_1(Y').$$ Under the cobordism $W$ induced
by the two-handle, the map $$F_{W}\colon \HFa(Y)\longrightarrow
\HFa(Y')$$
is specified by $$F_{W}(\xi\cm \Theta)=\pi(\xi)\cm
\Theta',$$
where here $\Theta'$ is some fixed generator of $\HFa(Y')$,
and $\xi$ is any element of $\Wedge^* H_1(Y)$.
Dually, if $K\subset Y$ is an unlink, then $Y''=Y_0(K)\cong
\#^{k+1}(S^2\times S^1)$, with a natural inclusion
$$i\colon H_1(Y) \longrightarrow H_1(Y'').$$
Under the cobordism $W'$ induced by the two-handle the map
$$F_{W'}\colon \HFa(Y) \longrightarrow \HFa(Y'')$$
is specified by
$$F_{W'}(\xi\cm \Theta)=\xi\wedge [K'']\cm \Theta'',$$
where here $[K'']\in H_1(Y'';\Z)$ is a generator in the kernel of 
the map $H_1(Y'')\longrightarrow H_1(W')$.
\end{prop}

\begin{proof}
  The identification of $\HFa(Y)$ follows from a direct inspection of
  the Heegaard diagram, as explained in
  Subsection~\ref{HolDiskTwo:subsec:GenusOne}.  The fact that
  $F_W(\Theta)$ is a generator for $\HFa(Y')$ (which we denote by
  $\Theta'$) follows from a direct inspection of a Heegaard triple
  which naturally splits into genus one summands, c.f.~\cite{HolDiskFour}.
  (Alternately, one could use the surgery exact sequence which in 
  this case reads
  $$
  \begin{CD}
    ...@>>>\HFa(Y') @>>>\HFa(Y) @>>>\HFa(Y')@>>>...
  \end{CD}
  $$
  to deduce that the map from $\HFa(Y)$ to $\HFa(Y')$, which in
  this case is induced by the two-handle $W$ equipped with its torsion
  $\SpinC$ structure, is surjective.)  The more general formula for
  $F_W$ follows from naturality of the triangle maps under the $H_1$
  action (c.f.~\cite{HolDiskFour}).  The case of $W'$ follows
  similarly.
\end{proof}

See~\cite{HolDisk} and~\cite{HolDiskFour}.

We can now link Khovanov's construction (using notation from
Section~\ref{sec:Khovanov}) with $\HFa$:

\begin{prop}
\label{prop:E2Term}
Fix a projection $\Diagram$ for $K$.
There is an isomorphism for each $I$
$$\Psi(I)\colon {\widetilde V}(\Diagram(I))\stackrel{\cong}\longrightarrow
\HFa(\Sigma(\Diagram(I))),$$
which is natural under cobordisms, in the following sense.
If $I'$ is an immediate successor of $I$, then there is a naturally
induced cobordism (induced from a single two-handle addition)
from $\Sigma(\Diagram(I))$ to $\Sigma(\Diagram(I'))$, and hence an
induced map
$$\Ga_{I<I'}\colon \HFa(\Sigma(I))\longrightarrow \HFa(\Sigma(I')).$$
Naturality of $\Psi$ is captured in the following commutative diagram,
which is valid whenever $I'$ is an immediate successor of $I$:
\begin{equation}
\label{eq:Naturality}
\begin{CD}
{\widetilde V}(\Diagram(I))@>{\dkh_{I<I'}}>> 
{\widetilde V}(\Diagram(I')) \\
@V{\psi(I)}VV @VV{\psi(I')}V \\
\HFa(\Sigma(\Diagram(I)))@>{G_{I<I'}}>> 
\HFa(\Sigma(\Diagram(I')))
\end{CD}
\end{equation}
\end{prop}

\begin{proof}
First, note that for each $I$, we can write
$\Diagram(I)=S_0\cup...\cup S_k$, where here the $S_i$ are pairwise
disjoint unknots, and $k\geq 0$.  In this case,
$\Sigma(\Diagram(I))\cong \#^{k}(S^2\times S^1)$.  Indeed, we give a
basis $\{[\gamma_i]\}_{i=1}^{k}$ for $H_1(\Sigma(\Diagram(I)))$ as
follows. For $i>0$, let $[\gamma_i]\in
H_1(\Sigma(\Diagram(I)))$ be the homology class of the curve obtained
as the branched double cover of an arc from $S_0$ to $S_i$ (recall
that we are using here $\Zmod{2}$ coefficients). This induces the
identification $${\widetilde Z}(\Diagram(I))\cong
H_1(\Sigma(\Diagram(I))).$$ Combined with
Proposition~\ref{prop:STwoSOne}, we get a canonical identification
$$\Vred(\Diagram(I))\cong \HFa(\Sigma(\Diagram(I))).$$ Commutativity
of Diagram~\eqref{eq:Naturality} is proved in four cases, each of
which follows from Proposition~\ref{prop:STwoSOne}.

Suppose that $I'$ is obtained from $I$ by merging two circles $S_1$
and $S_2$, neither of which is marked. Then, we claim that 
in the cobordism $W$, the curves $\gamma_1$ and $\gamma_2$ become homologous;
indeed, both are homologous to the new curve $\gamma_1'$. Commutativity
of the square now follows readily from Proposition~\ref{prop:STwoSOne} and
the definition of $\dkh_{I<I'}$.

Dually, when $I'$ is obtained from $I$ by splitting 
an unmarked circle $T_1$ into two circles $S_1$ and $S_2$,
the curve  $\gamma_1-\gamma_2$ is null-homologous in 
the induced cobordism $W'$. Again, commutativity of the claimed square
now follows readily from  Proposition~\ref{prop:STwoSOne}.

\begin{figure}
\mbox{\vbox{\epsfbox{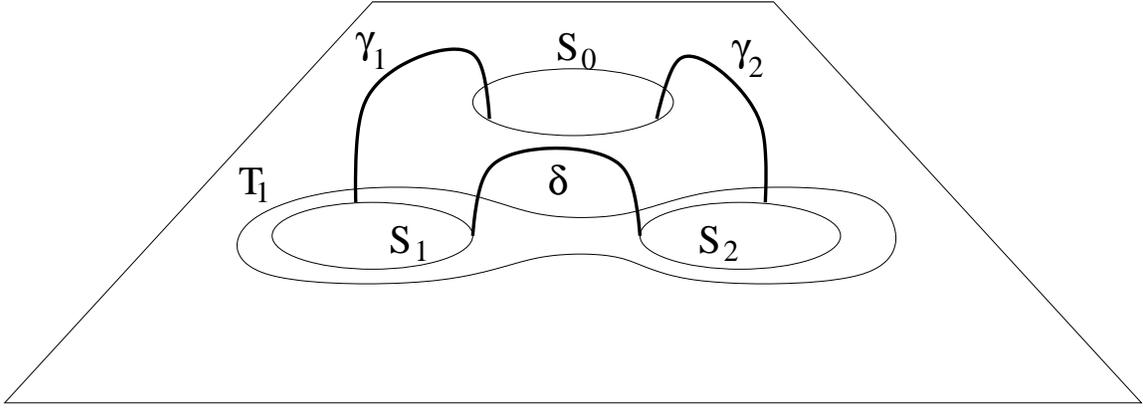}}}
\caption{\label{fig:TQFT}
  {\bf Homological relations in the cobordisms.}  The plane of the
  (un)-link projection is indicated by the quadrilateral, on which we
  nhave the marked component $S_0$, two unlink components $S_1$ and
  $S_2$, and an alternative component $T_1$, obtained by merging $S_1$
  and $S_2$. We have also illustrated the curves $\gamma_1$ and
  $\gamma_2$. The picture illustrates that $\gamma_1$ and $\gamma_2$
  become homologous after $S_1$ and $S_2$ are merged. Dually, it
  illustrates that in the cobordism where $T_1$ is divided in two, it
  is the curve $\delta=\gamma_1-\gamma_2$ which becomes
  null-homologous.}
\end{figure}

The two corresponding cases involving a marked circle follow similarly
(indeed, they follow formally in the same manner, once we declare
$[\gamma_0]=0$).
\end{proof}

With these preliminaries in place, we can now state and prove the
following precise version of Theorem~\ref{intro:DCovSpecSeq} in the
introduction.

\begin{theorem}
\label{thm:DCovSpecSeq}
Given a projection $\Diagram$ of a link $L$, there
is a spectral sequence converging to $\HFa(\Sigma(L))$
whose $(E^1,d^1)$ complex is 
isomorphic to 
Khovanov's reduced chain complex (for the mirror of $L$); 
i.e. there are isomorphisms
$\psi^m$ making the following diagram commute:
$$
\begin{CD}
E^1_{m-n_+(\Diagram)} @>{d^1_{m-n_+(\Diagram)}}>> E^1_{m-n_+(\Diagram)-1} \\
@V{\psi^m}VV @V{\psi^{m-1}}VV  \\
\CKhred(\Diagram(L),m) @>{\dkh^m}>>\CKhred(\Diagram(L),m-1).
\end{CD}
$$
In particular, the $E^2$ term of this sequence
is identified with Khovanov's homology of
$L$.
\end{theorem}

\begin{proof}
As explained in Section~\ref{sec:BranchedCovers}, a diagram $\Diagram$
for a link $K$ with $\ell$ crossings gives rise to a link $L$ in
$Y=\Sigma(K)$ whose components correspond to the crossings of
$\Diagram$. Moreover, for each $I\subset \{0,1\}^{\ell}$, the
three-manifold obtained by performing surgeries along these components
of $L$ is the branched double cover of $S^3$ branched along the
collection of unknots $\Diagram(I)$.  We now apply
Theorem~\ref{thm:LinkSpecSeq}.  The identification with Khovanov's
(reduced) complex is now provided by
Proposition~\ref{prop:E2Term}. Note that the spectral sequence coming
from Theorem~\ref{thm:LinkSpecSeq} is filtered by a cube, rather than
$\Z$.  We pass to a $\Z$-filtered object by ``flattening'' the cube as usual
(c.f. Remark~\ref{rmk:Flatten}).
\end{proof}

\vskip.3cm
\noindent{\bf{Proof of Corollary~\ref{cor:DetKhovBounds}.}}
The inequalities follow from 
$\det(L)=|H^2(\Sigma(L);\Z)|$ 
(c.f.~\cite{Lickorish}); but this agrees with $\chi(\HFa(\Sigma(L)))$
(c.f. Proposition~\ref{HolDiskTwo:prop:EulerHFa} of~\cite{HolDiskTwo}).
It follows at once that
$\det(L)\leq \rk \HFa(\Sigma(L))$.
The other inequality follows from Theorem~\ref{thm:DCovSpecSeq},
together with the straightforward inequality for spectral sequences
$\rk E^\infty \leq \rk E^2$.
\qed

\commentable{
\bibliographystyle{plain}
\bibliography{biblio}
}

\end{document}